\documentclass[12pt,thmsa,emstex]{article}
\UseRawInputEncoding
\usepackage{amsfonts}
\usepackage{t1enc}
\usepackage{amssymb}
\usepackage{epsfig}
\usepackage[french,english]{babel}
\usepackage{amsmath}
\usepackage{graphics}
\usepackage{layout}
\usepackage{latexsym}
\usepackage{color}
\usepackage{verbatim}
\def\clos{\overline}
\usepackage{xcolor} %oliver 
\long\def\metanote#1#2{{\color{#1}\ %oliver
\ifmmode\hbox\fi{\sffamily\mdseries\upshape [#2]}\ }} %oliver

\def\mylabel#1{\label{#1}}

\newtheorem{theorem}{Theorem}[section]
\newtheorem{lemma}[theorem]{Lemma}
\newtheorem{corollary}[theorem]{Corollary}
\newtheorem{proposition}[theorem]{Proposition}

\newtheorem{exercise}[theorem]{Exercise}
\newtheorem{remark}[theorem]{Remark}
\newtheorem{example}[theorem]{\bf{Example}}
\newtheorem{assumption}[theorem]{\bf{Assumption}}

\def\tcr{}
\def\tcb{}

\def\bit{\begin{itemize}}
\def\eit{\end{itemize}}

\reversemarginpar   %%puts label in LEFT margin

\def\bc{\begin{center}}
\def\ec{\end{center}}

\def\bthm{\begin{theorem}}
\def\ethm{\end{theorem}}

\def\bcor{\begin{corollary}}
\def\ecor{\end{corollary}}

\def\bprop{\begin{proposition}}
\def\eprop{\end{proposition}}

\def\blem{\begin{lemma}}
\def\elem{\end{lemma}}

\def\bex{\begin{example}}
\def\eex{\end{example}}

\def\bexo{\begin{exercise}}
\def\eexo{\end{exercise} }

\def\brem{\begin{remark}}
\def\erem{\end{remark}}
\def\prf{{\bf Proof: }}

\def\bdes{\begin{description}}
\def\edes{\end{description}}

\def\ita{\item[(a)]}
\def\itb{\item[(b)]}

\def\iti{\item[(i)]}
\def\itii{\item[(ii)]}

\def\itiii{\item[(iii)]}

\def\beq{\begin{equation}}
\def\eeq{\end{equation}}

\def\ben{\begin{enumerate}}
\def\een{\end{enumerate}}

\def\beqar{\begin{eqnarray}}
\def\eeqar{\end{eqnarray}}

\def\beqarr{\begin{eqnarray*}}
\def\eeqarr{\end{eqnarray*}}

\def\qed{\hfill $\Box$ \\[2ex]}
\def\prf{{\bf Proof: }\hspace{.1in}}

\catcode`\@=11

\newcommand{\C}{\mathcal{C}}
\newcommand{\M}{\mathcal{M}}

\newcommand{\EE}{\mathbb{E}}

\def\Ind{{\mathbf 1}}
\def\ZZ{{\mathbb Z}}       %bold Z
\def\RR{{\mathbb R}}  % bold R
\def\Rp{{\mathbb R}_+}   %bold R-plus

\def\NN{{\mathbb N}}

\def\Pr{{\mathsf P}}
\def\M{{\mathbf M}}
\def\m{{\mathbf m}}
\def\i{{\mathbf i}}

\def\rar{\rightarrow}
\def\eps{\varepsilon}

\begin{document}
\title{On invariant distributions of Feller Markov chains with applications to dynamical systems with random switching}
\author{Michel  Bena{\"i}m$^1$  and Oliver Tough$^2$}
\footnotetext[1]{Institut de Math\'{e}matiques, Universit\'{e} de Neuch\^{a}tel, Switzerland.}
\footnotetext[2]{Department of Mathematical Sciences, University of Bath, United Kingdom}
\maketitle

\begin{abstract}
We introduce simple conditions ensuring that invariant distributions of a Feller Markov chain on a compact Riemannian manifold are absolutely continuous with a lower semi-continuous, continuous or smooth density with respect to the Riemannian measure. This is applied to Markov chains obtained by random composition of maps and to piecewise deterministic Markov processes obtained by random switching between flows.
\end{abstract}
\tableofcontents
\section{Introduction}
 The aim of this paper is to propose and discuss simple conditions guaranteeing that the invariant distributions of a Feller Markov chain on a compact space satisfy certain regularity properties, such as having lower semi-continuous, continuous or smooth densities with respect to a reference measure.

Our initial motivation comes from {\em piecewise deterministic Markov processes} (PDMPs) generated by random switching between deterministic flows. The ergodic properties of these type of processes have been the focus of much attention in the last decade  and conditions ensuring existence, uniqueness, and absolute continuity (with respect to a reference Riemannian measure) of invariant measures, are now well understood (\cite{bakhtin&hurt}, \cite{BLMZ_2015},  \cite{BLMZ_2019}, \cite{Cloez15},  \cite{BCL17}, \cite{BHS18}). Concerning the regularity (continuity, smoothness) of these densities,  some partial results have been obtained in dimension one by Bakhtin, Hurth and Mattingly in \cite{bakhtin&hurt&matt}, Bakhtin, Hurth, Lawley and Mattingly in \cite{BHLM18} and \cite{BHLM21} for specific systems in dimension two, and by the present authors in \cite{BT23} for systems under ``sufficiently fast'' switching. Also worth mentioning is L\"ocherbach's beautiful article \cite{Loch18} on certain PDMPs with jumps, in which techniques (similar to those in \cite{BHLM18}) are used to prove regularity. However, beyond these cases,  the problem remains largely open. One of our principal goals is to revisit these questions, and to provide a simple and general framework allowing - in particular - for the results of \cite{BHLM18} and \cite{BT23} to be extended.

  The general  idea of the paper can be  roughly described as follows. Suppose $P$ is a Feller Markov kernel on some compact metric space $M$ and that ${\cal C}(M)$ is a convex cone of measures embedded in some Banach space $E.$  For instance, if $M$ is a Riemannian manifold, ${\cal C}(M)$ can be chosen to be the set of measures having a $C^r$ ($r \geq 0$) density with respect to the Riemannian measure, and  $E = C^r(M).$

  Suppose that $P = Q + \Delta$ where $Q, \Delta$ are sub-Markov kernels such that $Q$ maps the set of probability measures into ${\cal C}(M)$ and $\Delta$ maps ${\cal C}(M)$ into itself. Then, it is not hard to show that if the spectral radius of $\Delta$ (seen as an operator on $E$) is $< 1,$ invariant distributions of $P$ lie in ${\cal C}(M).$

 The paper explores and develops this idea. Section \ref{sec:notation} sets the general framework, notation and hypotheses. Here we state and prove our general results, such as the aforementioned Theorem \ref{th:main2}, along with other results ensuring absolute continuity of the invariant distributions and lower semi-continuity of their densities (Theorems \ref{th:main1} and \ref{th:piirreducible}).

 Section \ref{sec:RDS} considers the situation where $P$ is induced by a random iterative system on a compact Riemannian manifold and provides conditions ensuring that the decomposition $P = Q + \Delta$ holds with ${\cal C}(M)$ the set of measures having a density (respectively a lower semi-continuous, of $C^r$ density) with respect to the Riemannian measure. In the specific case where $\Delta = \delta_{\phi}$ with $\phi$ a local diffeomorphism, $\Delta$ is nothing but the Ruelle transfer operator of $\phi$ and its spectral radius can be estimated in terms of certain topological (or measure-theoretic) invariants for $\phi.$ This is done in Subsection \ref{sec:diffeo} and applied to specific examples in Subsection \ref{sec:applRDS}.

 Section \ref{sec:pdmp} is devoted to PDMPs, as described above. We prove that under certain Hörmander conditions, there are finitely many ergodic measures that are absolutely continuous with respect to the Riemannian measure and whose densities are  lower semi-continuous (Theorem  \ref{th:lscpdmp1}). If the Hormander condition holds at an accessible point, such a measure is unique (Theorem  \ref{th:lscpdmp0}). In Subsection \ref{sec:torus} we consider the situation of two transverse vector fields on the torus, and give a precise condition (involving the switching rates and the Floquet exponents of the linearly stable periodic orbits of the vector fields) ensuring that the invariant measures have a $C^k$ density (Theorem \ref{th:torus}). This result  relies on the spectral radius estimate of the Ruelle transfer operator given in section \ref{sec:diffeo} and substantially extends the results in \cite{BHLM18}. The last section \ref{sec:fast} is devoted to general PDMPs under fast switching. We show how our approach provides for a short proof that under fast switching and a certain Hörmander condition, invariant densities are $C^r$.

\section{Notation, hypotheses and basic results}
\label{sec:notation}
Let $M$ be a compact metric space equipped with its Borel sigma field ${\cal B}(M).$

We let ${\cal M}(M)$ (respectively  ${\cal P}(M)$) denote the set of non negative finite (respectively, probability) measures over $M.$

A {\em convex cone} of a measures is a set $\C(M) \subset {\cal M}(M)$  such that  $\alpha \mu + \beta \nu \in \C(M),$ for all $\mu, \nu  \in \C(M),$ and all $\alpha, \beta  \geq 0.$
 \bex
 \label{ex:riemann}
 {\rm Suppose that $M$ is a Riemannian manifold with Riemannian  measure $m.$ Examples of convex cones in ${\cal M}(M)$ include:

 \begin{itemize}
   \item ${\cal M}_{ac}(M) \subset {\cal M}(M),$ the set of measures which are absolutely continuous with respect to $m;$
   \item ${\cal M}_{ac}^{ls}(M) \subset {\cal M}_{ac}(M)$, the subset which have a lower semi-continuous density;
   \item ${\cal M}_{ac}^{r}(M) \subset {\cal M}_{ac}(M), r \geq 0,$ the subset which have a $C^r$ density.
 \end{itemize}
 }
 \eex
 \tcr{Here and throughout, when we say that a measure has a $l.s.c$ or $C^r$ density, we actually mean that a version of the density can be chosen to be  $l.s.c$ or  $C^r.$}

 A  {\em bounded kernel}  on $M$ is a  family $Q =\{Q(x,\cdot)\}_{x \in M}$ with $Q(x,\cdot) \in {\cal M}(M)$ such that for all $A \in {\cal B}(M)$, the mapping $x \rightarrow Q(x,A)$ is measurable, and  $\sup_{x \in M } Q(x,M) < \infty.$
 We say that $Q$ is  {\em non-degenerate} if $Q(x,M) > 0$ for all $x \in M;$   {\em sub-Markov} if $\sup_{x \in M } Q(x,M) \leq 1;$ and {\em Markov} if  $Q(x,\cdot) \in {\cal P}(M)$ for all $x \in M.$

We let  $B(M)$ (respectively $C^0(M)$) denote the Banach
space of bounded measurable (respectively continuous) real valued   functions on $M,$ endowed with the uniform norm
$\|f\|_{0} = \sup_{x \in M} |f(x)|.$

A   bounded kernel $Q$ induces a bounded operator on $B(M)$ defined  by $$Qf(x) = \int_M f(y) Q(x,dy),$$ for all $f \in B(M).$ We call it  {\em Feller} if it  maps $C^0(M)$ into itself.
It also induces an operator on ${\cal M}(M)$  defined by $$\mu Q(A) = \int \mu(dx) Q(x,A),$$ for all $\mu \in {\cal M}(M)$ and
$A \in {\cal B}(M).$

  If $Q$ is   Markov, we let
 $\mathsf{Inv}(Q)$ denote  the set of {\em invariant probability measures} of $Q.$; that is the set of $\mu  \in {\cal P}(M)$  such that  $\mu Q= \mu.$ If  $Q$ is Markov and Feller, then
  $\mathsf{Inv}(Q)$ is a non-empty convex  compact (for the weak* topology) subset of ${\cal P}(M)$ (see e.g \cite{BenaimHurth}, Corollary 4.21).
%%%%%%%%%%%%%%%%
\medskip

From now on, we let $P$ denote a Markov Feller kernel and $\C(M)$ a convex cone of measures.
Our standing assumption  is given by the following assumption.
  \begin{assumption}[Standing assumption]
  \label{hyp1} The kernel $P$ may be decomposed into $P = Q + \Delta,$ \tcb{where}:
  \bdes
  \iti $Q$ is a  non-degenerate  Feller sub-Markov kernel and $\Delta$ is a (possibly degenerate) sub-Markov kernel;
   \itii ${\cal M}(M) Q := \{  \mu Q: \: \mu \in {\cal M}(M) \} \subset \C(M);$
   \itiii $\C(M) \Delta := \{  \mu \Delta: \: \mu \in \C(M) \} \subset \C(M).$
   \edes
  \end{assumption}
   In our applications, $\C(M)$ will be,  like in Example \ref{ex:riemann},  a set of measures having certain regularity properties. In words, Assumption \ref{hyp1} means that $Q$ "creates" regularity, whilst $\Delta$ "preserves" regularity.

\tcr{Before going further, it is worth pointing out that  the idea to decompose $P$ as $P = Q + \Delta$ where $Q$ enjoys certain regularity properties is part of the folklore in the literature on Markov chains. It is reminiscent of the minorization condition (in this case $Q(x,\cdot) := \nu(\cdot)$) introduced in the late 70's by Athreya and Ney \cite{AN78} in their analysis of Harris chains (see also  Meyn and Tweedie \cite{MT93}, or Duflo \cite{duf00}). In case  $Q$ is a  {\em continuous component} (meaning that $x \mapsto Q(x,A)$ is lower semi-continuous for all Borel set $A$) we retrieve the notion of {\em $T$-chain} introduced by  Meyn and Tweedie \cite{MT93}, Chapter 6.}

It follows from Assumption \ref{hyp1} that $\Delta$ is  Feller and that $$\sup_{x \in M} \Delta(x,M) := \rho < 1.$$
In particular,  $$(I-\Delta)^{-1} :=  \sum_{k \geq 0} \Delta^k$$ is also a  Feller kernel and $$(I-\Delta)^{-1}(x,M) \leq \frac{1}{1-\rho}.$$ Here $I = \Delta^0 = \{\delta_x(\cdot)\}_{x \in M}.$

The following result is a straightforward consequence of Assumption \ref{hyp1}, and will be used repeatedly.
\blem
\label{lem:pi}
Let $\Pi \in \mathsf{Inv}(P).$ Then, under  Assumption \ref{hyp1} $(i),$ $$\Pi = \Pi Q (I- \Delta)^{-1} = \sum_{k \geq 0} \Pi Q \Delta^k.$$
\elem
\prf This follows directly from the equation $\Pi = \Pi P \Leftrightarrow \Pi (I - \Delta) = \Pi Q.$ \qed
 \bex
 \label{ex:sampling}
{\rm
Suppose that $Q(x,dy) = \pi(dy)$ with $\pi \in {\cal M}(M).$ Lemma \ref{lem:pi} shows that $$\mathsf{Inv}(P) = \{ \pi (I-\Delta)^{-1}\}.$$}
\eex

 We say that $\C(M)$ is {\em stable by monotone convergence} if for every sequence $(\mu_n)_{n \geq 0}$ with $\mu_n \in \C(M)$ and  $\mu_n \leq \mu_{n+1},$ $\mu := \lim_{n \rar \infty} \mu_n$ lies in $\C(M).$ Here, $\mu := \lim_{n \rar \infty} \mu_n$  simply means that  $\mu(A) := \lim_{n \rar \infty} \mu_n(A) \in [0,\infty]$ for all $A \in {\cal B}(M).$

\brem
\label{rem:monot}
{\rm The sets ${\cal M}_{ac}(M)$ and  ${\cal M}_{ac}^{ls}(M)$ as defined in Example \ref{ex:riemann} are stable by monotone convergence.
This is obvious for the first set. For the second, let $f_n$ be an l.s.c version of the density $\mu_n.$ Set $g_n= \max(f_1, \ldots ,f_n).$ Then $g_n$ is l.s.c and a density of $\mu_n$, using that $(\mu_n)$ is monotone so $g_n=f_n$ almost surely. Clearly also $(g_n)$ is monotone, and the limit is an l.s.c version of $\mu.$}
\erem
A first useful (and immediate) consequence of Lemma \ref{lem:pi} is the next result.
  \bthm
  \label{th:main1}
  Assume Assumption \ref{hyp1} holds with  $\C(M)$  stable by  monotone convergence. Then $\mathsf{Inv}(P) \subset \C(M).$
  \ethm

\bcor
\label{cor:LSC}
Suppose $M$ is a Riemannian manifold.   Assume Assumption \ref{hyp1} holds with $\C(M) =  {\cal M}_{ac}^{ls}(M).$ Then:
\bdes
\iti $\mathsf{Inv}(P) \subset \C(M);$
 \itii if $\mu, \nu \in \mathsf{Inv}(P)$ are  ergodic,  either $\mu = \nu$ or there exist nonempty disjoint  open sets $U, V$ such that $\mu(U) = \nu(V) = 1.$
In particular, if $M$ is connected and an invariant distribution has full support, then it is the unique invariant distribution of $P.$
\edes
\ecor
\prf $(i)$ follows from Proposition \ref{th:main1} and Remark \ref{rem:monot}. We now turn to $(ii).$ By ergodicity either $\mu = \nu$ or $\mu$ and $\nu$ are mutually singular. By Proposition  \ref{th:main1}, $\mu(dx) = h(x) m(dx)$ and $\nu(dx) = g(x) m(dx)$ with $h$ and $g$ lower semi-continuous. Set $U := \{x \in M \: : h(x) > 0\}$ and $V := \{x \in M \: : g(x) > 0\}.$ Then $U$ and $V$ are open and  $\mu(dx) \geq \frac{h(x)}{g(x)} 1_{V}(x) \nu(dx).$ Thus, if $\mu$ and $\nu$ are mutually singular, $h$ has to be zero on $V.$
\qed

Another useful (and immediate) consequence of Lemma \ref{lem:pi} is given by the next result.
\bthm
\label{th:main2}
Assume Assumption \ref{hyp1} holds with $\C(M)$  a closed subset of some  Banach space $(E, \|.\|_E).$ Assume furthermore  that the two following conditions hold:
\bdes
\iti $\sum_{k \geq 0} \|\mu \Delta^k\|_E < \infty$ for all  $\mu \in \C(M)$,
\itii For every Borel set $A \subset M,$ the map $\C(M) \rar \RR,  \mu  \mapsto \mu(A),$ is continuous when $\C(M)$ is equipped with the distance induced by $\|\cdot\|_E.$
\edes

 Then $\mathsf{Inv}(P) \subset \C(M).$
\ethm
\brem
\label{rem:Banach}
{\rm In the following sections, this theorem will be used when $M$ is a Riemannian manifold, $\C(M) = {\cal M}_{ac}^r(M)$, and $E$ is the Banach space of bounded signed measures whose density is $C^r$ (naturally identified with $C^r(M)$ equipped with the $C^r$ norm).}
\erem
\brem
{\rm A sufficient practical condition ensuring condition $(i)$  in Theorem \ref{th:main2} is that  $\mu \rar \mu \Delta$ extends to a bounded operator on $E$ whose spectral radius, $${\cal R}(\Delta,E) = \lim_{n \rar \infty} \|\Delta^n\|_E^{1/n},$$ is  strictly less than $1.$
}
\erem
\subsection{On Assumption \ref{hyp1}: a uniqueness result}
\label{sec:onH1H2}
It is often the case that a Markov kernel $P$ doesn't satisfy the standing assumption, Assumption \ref{hyp1}, but that some power of $P, P^k$ (for some $k \geq 1$), or its $a$-resolvent
$$R_a = (1-a) \sum_{k \geq 0} a^k P^k$$ (for some $0 < a < 1$), does. Since $$\mathsf{Inv}(R_a) = \mathsf{Inv}(P) \subset \mathsf{Inv}(P^k),$$ the conclusions of the previous theorems  remain valid in these cases.

 The next theorem illustrates this idea. Let  $P$ be a Feller Markov kernel which doesn't necessarily satisfy the standing assumption.
 A point $p \in M$ is called {\em accessible} (for $P$) if for every neighbourhood $U$ of $p$ and every $x \in M,$ $R_a(x,U) > 0$ (for some, hence all $0 < a < 1$).  The set of points which are accessible for $P$ is then the (possibly empty) compact set
 $$\Gamma_P = \bigcap_{x  \in M } \mathsf{supp}(R_a(x,\cdot)),$$
 \tcr{where $\mathsf{supp}(R_a(x,\cdot))$ stands for the topological support of the measure $R_a(x,\cdot).$}
Point $p$ is  called a {\em weak Doeblin} point if there exists a neighbourhood $V$ of $p,$  a non-trivial measure
$\pi \in {\cal M}(M)$, and $0 < a < 1,$ such that $ R_a(x,dy) \geq \pi(dy)$ for all $x \in V.$ The measure $\pi$ is called a
{\em minorizing} measure.
\bthm
 \label{th:piirreducible} Let $\C(M)$ be  a convex cone stable by monotone convergence. Suppose that $\C(M) P \subset \C(M)$ and that $P$ possesses  an accessible weak Doeblin point with a minorizing measure $\pi \in \C(M)$ such  that for all $\mu \in \C(M)$, $$\mu \geq \pi \Rightarrow \mu - \pi \in \C(M).$$  Then $P$ has a unique invariant probability measure $\Pi$ and
 $\Pi \in \C(M).$
\ethm
\prf By assumption, there exists an open set $V$ such that $R_a(x,dy) \geq \pi(dy)$ for all $x \in V$ and $R_a(x,V) > 0$ for all $x \in M.$ By the Feller continuity of $P$ (hence of $R_{a}$), \tcr{$x \rar R_a(x,V)$ is lower semi-continuous. Then, by  compactness}  $R_{a}(x,V) \geq \delta > 0$ for all $x \in M$ and some $\delta > 0.$ It follows that $R_a^2(x,dy) \geq \delta \pi(dy).$ By Theorem \ref{th:main1} and Example \ref{ex:sampling} applied to $R_a^2$, with $Q = \delta \pi, \Delta = R_a^2 - \delta \pi,$ we get that $\mathsf{Inv}(P) \subset \mathsf{Inv}(R_a^2) = \{\Pi\} \subset \C(M).$ \qed
\tcr{Note that the minoration $R_a^2(x,dy) \geq \delta \pi(dy)$ in the proof above, implies that $P$ is {\em $\psi$-irreducible} in the sense of Meyn and Tweedie \cite{MT93}, and it is well known that a $\psi$-irreducible chain has (at most) one invariant probability measure (see e.g~\cite{MT93}, \cite{duf00}, \cite{BenaimHurth}). The added value of Theorem \ref{th:piirreducible} is the simple proof that  $\Pi \in \C(M).$}
\section{Random maps}
\label{sec:RDS}
We suppose here that $M$ is a compact $d$-dimensional connected Riemannian  manifold. For $k \geq 0,$ we let $C^k(M)$ denote the space of $C^k$ functions $\rho : M \rar \RR,$ equipped with the $C^k$  topology (see e.g~\cite{Hirsch76}, Chapter 2).
We let $\|\cdot\|_{C^k(M)}$ denote a norm on $C^k(M)$ making $C^k(M)$ a Banach space. We let $C^k(M,M)$ be the space of $C^k$ maps from $M$ into itself, equipped with the $C^k$ topology and associated Borel $\sigma$-field.
%Without loss of generality we can always assume that $M$ is a sub-manifold of $\RR^N$ for $N$ sufficiently large and see $C^k(M,M)$ as a subset of $C^k(M, \RR^N) = C^k(M)^N.$

We now let $r \geq 1,$ and let $\nu$ be a probability measure on $C^r(M,M).$ Consider the chain on $M$ induced by the random iterative system
$$X_{k+1} = \varphi_{k+1}(X_k),$$ where $(\varphi_k)_{k \in \NN}$ is a family of i.i.d random variables, independent of $X_0,$  having distribution $\nu.$
%$\varphi : M \times   \RR^n \mapsto M, (x,\theta) \rar \varphi(x,\theta) = \varphi_{\theta}(x)$ be a $C^r, r \geq 1,$ map.

The kernel of this chain can then be written
\beq
\label{eq:kernelRI} P^{\nu}f(x) = \int_{C^r(M,M)} f(\varphi(x)) \nu(d\varphi),
 \eeq and is clearly Feller. For further reference, we call this kernel the {\em kernel induced by $\nu.$ }

Throughout this section we shall take $P := P^{\nu}$, and assume that $\nu$ may be written as
$$\nu := (1-a) \nu_0 + a \nu_1,$$ \tcb{where}  $\nu_0, \nu_1$ are two probability measures over $C^r(M,M)$ and $0 < a < 1.$ Thus we can write $P = Q + \Delta$ with $$Q = (1-a) P^{\nu_0}$$ and $$\Delta = a P^{\nu_1},$$ where $P^{\nu_0}, P^{\nu_1}$ are defined like $P^{\nu}$ with $\nu_0, \nu_1$ in place of $\nu.$  We furthermore assume that $\nu_0, \nu_1$ satisfy the following hypotheses \ref{hyp:RDS1} and \ref{hyp:RDS2} below. These are natural hypotheses ensuring that the standing assumption, Assumption \ref{hyp1}, holds true with $\C(M)$ being one of the sets ${\cal M}_{ac}(M), {\cal M}_{ac}^{ls}(M)$ or ${\cal M}_{ac}^{r-1}(M)$ as defined in Example \ref{ex:riemann}. To be concise, Assumption \ref{hyp:RDS1} assumes that $\nu_0$ is the image measure of a finite dimensional $C^r$ density by a submersion, while Assumption \ref{hyp:RDS2} assumes that $\nu_1$ is supported by local diffeomorphisms.
\begin{assumption}[Standing assumption $1$ for RDS]
\label{hyp:RDS1}
There exist $n \geq d,$  a smooth $n$-dimensional manifold $\Theta$ with smooth Riemann measure $d\theta,$ a $C^r$ probability density function $h_0 : \Theta \rar \Rp$ with compact support $\mathsf{supp}(h_0),$ and a $C^r$ map
$${\mathbf{\Phi}}: M \times \Theta \mapsto M,$$ $$(x,\theta)  \rar \mathbf{\Phi}(x,\theta) = \mathbf{\Phi}_{\theta}(x)$$ such that:
\bdes
\iti $\nu_0$ is the image measure of $h_0(\theta) d\theta$ by the map $\theta \rar \mathbf{\Phi}_{\theta}.$ That is
$$P^{\nu_0}(f)(x) = \int_{\Theta} f(\mathbf{\Phi}_{\theta}(x)) h_0(\theta) d\theta.$$
\itii $\partial_{\theta} \mathbf{\Phi}(x,\theta) : T_{\theta} \Theta \mapsto T_{\mathbf{\Phi}(x,\theta)} M$ is surjective for all $x \in M$ and $\theta \in \mathsf{supp}(h_0)$.
\edes
\end{assumption}

\tcr{The next proposition relies on the  fact that the push-forward of a measure having a smooth, compactly supported density by a smooth submersion has a smooth density. It is reminiscent of Lemma $6.3$ in  \cite{BLMZ_2015} and Lemma $2$ in \cite{bakhtin&hurt}.}
\bprop
\label{lem:submersivemap} Assume Assumption \ref{hyp:RDS1}. Then,
there exists a $C^r$ map $q : M \times M \mapsto \Rp$  such that  $$P^{\nu_0}(x,dy) = q(x,y) m(dy).$$
In particular, ${\cal M}(M)Q \subset  {\cal M}_{ac}^{r}(M).$
\eprop
\prf We assume for notational convenience that $\Theta = \RR^n$, but the proof easily extends to the general case.

\underline{Claim}: For all $x^* \in M$ and  $\theta^* \in \mathsf{supp}(h_0),$ there exist neighbourhoods $U (= U(x^*,\theta^*))$ of $x^*$ and  $V (= V(x^*,\theta^*))$ of $\theta^*$ such that for every $C^r$ function $\eta : \RR^n \mapsto \RR$ with compact support  $\mathsf{supp}(\eta) \subset V,$ there exists a $C^r$ map $q_{\eta} : M \times M \rar \Rp$ with the property that
  $$\int_{\RR^n} f(\mathbf{\Phi}(x,\theta)) h_0(\theta) \eta(\theta)  d\theta = \int_M q_{\eta}(x,y)f(y)m (dy)$$ for all $x \in U,$ and $f \in B(M).$

We assume for the time being that the claim is proven. Fix $x^* \in M.$ We extract from the family $\{V(x^*,\theta^*), \theta^* \in \mathsf{supp}(h_0)\}$ a covering of $\mathsf{supp}(h_0)$ by open sets $V_i = V(x^*, \theta_i), i \in I$, with $I$ finite. Set $U = \cap_{i \in I} U(x^*, \theta_i).$ Using a partition of unity subordinate to $\{V_i\}_{i \in I},$ $h_0$ can be written as
 $h_0  = \sum_{i \in I}  h_0 \eta_i$ where $\eta_i$ is smooth with compact support in $V_i,$ $0 \leq \eta_i$, and $\sum_{i \in I} \eta_i = 1.$ It then follows from the claim that for all $x \in U,$
 $$P^{\nu_0}(x,dy) = \sum_{i \in I} q_i(x,y) m(dy),$$  where $q_i : M \times M \rar \Rp$ is $C^r.$ This proves the proposition.

   \underline{Proof of the claim}: After a permutation of the canonical basis of $\RR^n$ we can assume that $\theta =  (\theta_1,\theta_2) \in \RR^d \times \RR^{n-d}$ where $\partial_{\theta_1} \mathbf{\Phi}(x^*,\theta^*)$ has rank $d.$  Thus, by the inverse function theorem, there exist open neighbourhoods $U'$ of $x^*$ and $V = V_1 \times V_2$ of $\theta^* = (\theta_1^*,\theta_2^*)$
    such that the map $$H: (\theta_1,\theta_2,x) \rar (\mathbf{\Phi}(x,\theta), \theta_2, x)$$ is a $C^r$ diffeomorphism from $V \times U'$ onto its image $W = H(V \times U').$ Its inverse is then given by $(y,\theta_2,x) \rar (\psi(y,\theta_2,x), \theta_2,x)$, where $\psi : W  \mapsto V_1$ is $C^r.$

   Let $U$ be a  neighbourhood of $x^*$ with $\overline{U} \subset U'$, and let $\eta : \RR^n \rar \Rp$ be a $C^r$ function with compact support $\mathsf{supp}(\eta) \subset V.$ Set $K = \mathsf{supp}(\eta) \times \overline{U}$ and let $\tilde{k}(x,y,\theta_2)$ be a $C^r$ function which coincides with $$(\eta h_0)(\psi(y,\theta_2,x),\theta_2)|\mathsf{det} \partial_{y} \psi(y,\theta_2,x)|$$ on $H(K)$ and is zero outside $W.$
 We define $q_{\eta} : M \times M \mapsto \Rp$ by $$q_{\eta}(x,y) = \int \tilde{k}(x,y,\theta_2) d\theta_2.$$  Then $q_{\eta}$ is $C^r$ and by the change of variable formula,
   $$\int f(\mathbf{\Phi}(x,\theta)) (\eta h_0)(\theta) g(x) d\theta m(dx) = \int q_{\eta}(x,y) g(x) f(y) m(dx) m(dy)$$
   for every continuous function $g$ with support contained in $U.$ This proves the claim. \qed
   %%%%%%%%%%%

    We define $\mathsf{Diff_{loc}^r}(M) \subset C^r(M,M)$ to be the (open) set of maps $\varphi \in  C^r(M,M)$ for which $D\varphi(x) : T_xM \mapsto T_{\varphi(x)}M$ is invertible at every point  $x \in M.$

We let $\varphi \in \mathsf{Diff_{loc}^r}(M).$ It is not hard to see that  $\varphi^{-1}(y)$ is nonempty, finite, and that its cardinality doesn't depend on $y$ for all $y \in M$. \tcr{Indeed, by the  inverse function theorem, for each $x \in \varphi^{-1}(y), \varphi$ is a diffeomorphism from a neighborhood of $x$ onto a neighborhood of $y.$ This makes  $\varphi^{-1}(y)$ finite (by compactness) and  the mapping $y \rar \mathsf{card}(\varphi^{-1}(y))$  locally constant. By connectedness, it is constant.} We denote this cardinality by $\mathsf{deg}(\varphi)$.

We let $J(\varphi,x) > 0$ denote the {\em Jacobian} of $\varphi$ at $x$ with respect to $m.$ If the tangent spaces $T_x M$ and $T_{\varphi(x)}M$ are equipped with orthonormal bases, then  $$J(\varphi,x) = |\mathsf{det} D\varphi(x)|.$$
  The {\em transfer} or {\em Ruelle-Perron-Frobenius} operator induced by $\varphi$ is the operator ${\cal L}_{\varphi}$ acting on $L^1(m)$ or $C^{r-1}(M)$, defined by
  \beq
   \label{eq:transfer}
   {\cal L}_{\varphi}(\rho)(y) = \sum_{\{x \in \varphi^{-1}(y)\}} \frac{\rho(x)}{J(\varphi, x)}.
   \eeq
   \tcr{This definition is motivated by the change of variable formula. Indeed, if a measure has density $\rho$, its image-measure by $\varphi$ has density  ${\cal L}_{\varphi}(\rho).$}
   The fact that ${\cal L}_{\varphi}(\rho)$ maps  $C^{r-1}(M)$ into itself easily follows from the inverse function theorem. Indeed, for all $y \in M,$ there exist an open neighbourhood $U$ of $y$ and $C^r$ diffeomorphisms $\psi_i : U \mapsto \psi_i(U), i = 1, \ldots, \mathsf{deg}(\varphi)$, such that for all $z \in U,$
   $$ {\cal L}_{\varphi}(\rho)(z) = \sum_{i =1}^{\mathsf{deg}(\varphi)} \frac{\rho(\psi_i(z))}{J(\varphi, \psi_i(z))}.$$
   This expression also shows that  ${\cal L}_{\varphi}$ is a bounded operator on $C^{r-1}(M).$ We let $$\|{\cal L}_{\varphi}\|_{C^{r-1}(M)} = \sup_{\{\rho \: : \|\rho\|_{C^{r-1}(M)} \leq 1\}} \|{\cal L}_{\varphi}(\rho)\|_{C^{r-1}(M)}$$ denote its operator norm.

   For $0 \leq k \leq r-1,$ we let
   \beq
   \label{eq:RLphi}
   {\cal R}({\cal L}_{\phi}, C^k(M)) = \lim_{n \rar \infty} \|({\cal L}_{\phi})^n\|_{C^k(M)}^{1/n}
   \eeq be the spectral radius of ${\cal L}_{\phi}$ on $C^k(M).$
  \begin{assumption}[Standing assumption $2$ for RDS]
  \label{hyp:RDS2}
  $$\nu_1(\mathsf{Diff_{loc}^r}(M)) = 1.$$
  %  $$ \int_{\mathsf{Diff_{loc}^r}(M)} \|{\cal L}_{\varphi}\|_{C^{r-1}(M)} \nu_1(d\varphi) < \infty.$$
  \end{assumption}
\bprop
\label{lem:Lnu1} Assume Assumption \ref{hyp:RDS2}.
If $\mu \in {\cal M}_{ac}(M)$ has density $\rho,$ then  $\mu P^{\nu_1} \in {\cal M}_{ac}(M)$ and its density is given by
   $$y \rar {\cal L}_{\nu_1}(\rho)(y): =  \int_{\mathsf{Diff_{loc}^r}(M)} ({\cal L}_{\varphi} \rho)(y) \nu_1(d\varphi).$$
   This density  is lower semi-continuous whenever $\rho$ is. In particular $\C(M) \Delta \subset \C(M)$ with $\C(M) = {\cal M}_{ac}^{ls}(M)$ \tcr{where we recall (see the beginning of Section \ref{sec:applRDS}) that $\Delta = a P^{\nu_1}$ .}

   If in addition $$ \int_{\mathsf{Diff_{loc}^r}(M)} \|{\cal L}_{\varphi}\|_{C^{r-1}(M)} \nu_1(d\varphi) < \infty,$$ then ${\cal L}_{\nu_1}$ is a bounded operator on $C^{r-1}(M)$ and
   $$\|{\cal L}_{\nu_1}\|_{C^{r-1}(M)} \leq \int_{\mathsf{Diff_{loc}^r}(M)} \|{\cal L}_{\varphi}\|_{C^{r-1}(M)} \nu_1(d\varphi).$$
   In particular $\C(M) \Delta \subset \C(M)$ with  $\C(M) = {\cal M}_{ac}^{r-1}(M).$
\eprop
\prf For all $f \in B(M),$
\begin{align*}
\int_M P^{\nu_1}(f)(x) \rho(x) m(dx)  &=  \int_{M} \left( \int_{\mathsf{Diff_{loc}^r}(M)} f (\varphi(x)) \nu_1(d\varphi) \right) \rho(x) m(dx)\\
&= \int_{\mathsf{Diff_{loc}^r}(M)} \left( \int_{M} f (\varphi(x)) \rho(x) m(dx)\right )  \nu_1(d \varphi)
 \\
 &=  \int_{\mathsf{Diff_{loc}^r}(M)}  \left(\int_{M}  f(x) {\cal L}_{\varphi}(\rho)(x) m(dx)\right)  \nu_1(d \varphi)\\
 &=\int_M f(x) \left (\int_{\mathsf{Diff_{loc}^r}(M)} ({\cal L}_{\varphi} \rho)(x) \nu_1(d\varphi)\right) \rho(x) m(dx).
\end{align*}

The second and last equalities follow from Fubini's theorem, and the third one follows from the change of variable formula. This proves the first assertion.
 If $\rho$ is lower semi-continuous, so is ${\cal L}_{\varphi} \rho.$ Thus, if $y_n \rar y,$
  $$\liminf_{n \rar \infty} \int {\cal L}_{\varphi} \rho(y_n) \nu_1(d\varphi) \geq \int \liminf_{n \rar \infty} {\cal L}_{\varphi} \rho(y_n) \nu_1(d\varphi) \geq \int {\cal L}_{\varphi} \rho(y) \nu_1(d\varphi)$$ by Fatou's Lemma. This shows that  $\frac{d \mu P^{\nu_1}}{d m}$  is lower-semicontinuous.

  We now prove the last statement. For all $\rho \in C^{r-1}(M),$ the mapping ${\cal L}_{(\cdot)} \rho : \mathsf{Diff_{loc}^r}(M) \rar  C^{r-1}(M), \varphi \mapsto {\cal L}_{\varphi} \rho$ is continuous, hence measurable. It is then Bochner measurable (see   \cite{diestel}, Theorem 2, Section 1, Chapter 2) and the condition that $ \int_{\mathsf{Diff_{loc}^r}(M)} \|{\cal L}_{\varphi}(\rho)\|_{C^{r-1}(M)} \nu_1(d\varphi) < \infty$ makes it Bochner integrable (\cite{diestel}, Theorem 2, Section 2, Chapter 2). Properties of Bochner integrals  (\cite{diestel}, Theorem 4, Section 2, Chapter 2) imply that $$ \Big\|\int_{\mathsf{Diff_{loc}^r}(M)} {\cal L}_{\varphi}(\rho) \nu_1(d\varphi)\Big \|_{C^{r-1}(M)} \leq \int_{\mathsf{Diff_{loc}^r}(M)} \|{\cal L}_{\varphi}(\rho)\|_{C^{r-1}(M)} \nu_1(d\varphi).$$ This concludes the proof. \qed

We recall that $P = P^{\nu}$ is given by (\ref{eq:kernelRI}). Corollary \ref{cor:LSC} and Theorem \ref{th:main2} applied to the present setting, combined with Propositions \ref{lem:submersivemap} and \ref{lem:Lnu1}, imply the following.
\bthm
\label{th;main1RDS} Assume Hypotheses \ref{hyp:RDS1} and \ref{hyp:RDS2}. Then
$\mathsf{Inv}(P) \subset  {\cal M}_{ac}^{ls}(M).$  If $\mu \in \mathsf{Inv}(P)$ has full support, then $\mathsf{Inv}(P) =  \{\mu\}.$
\ethm
\tcr{We recall from the beginning of Section \ref{sec:applRDS} that $\nu = (1-a)\nu_0 + a \nu_1.$}
\bthm
\label{th:main2RDS} Assume Hypotheses \ref{hyp:RDS1} and \ref{hyp:RDS2}.
If  $$ \int_{\mathsf{Diff_{loc}^r}(M)} \|{\cal L}_{\varphi}\|_{C^{r-1}(M)} \nu_1(d\varphi) < \infty,$$
and $1/a$ is in the resolvent set of ${\cal L}_{\nu_1}$ (on $C^{r-1}(M)$), then $\mathsf{Inv}(P) \subset  {\cal M}_{ac}^{r-1}(M).$
\ethm
\subsection{Expansion volume rates and spectral radius}
%When $\nu_1$ is a Dirac}
\label{sec:diffeo}
In this subsection and the following, we consider the case where  $$\nu_1 = \delta_{\phi}$$ for some $\phi \in \mathsf{Diff_{loc}^r}(M), r \geq 1,$ so that ${\cal L}_{\nu_1}$ is the transfer operator ${\cal L}_{\phi}.$
 When $\phi$ is an {\em expanding} map (see the definition below), the spectral properties of ${\cal L}_{\phi}$ have been well understood since the seminal work of Ruelle \cite{Ruelle89}. We refer the reader to the excellent monograph \cite{Baladi2000} for a comprehensive introduction to the subject.

 When $\phi$ is non-expanding, it is still possible to give simple sufficient conditions ensuring that $\frac{1}{a}$ lies in the resolvent of ${\cal L}_{\phi},$ so that Theorem \ref{th:main2RDS} applies. This is the object of the next proposition, Proposition \ref{prop:transfert}.  Before stating this proposition we introduce certain quantities that will naturally appear in the estimate of the spectral radius of ${\cal L}_{\phi}:$ the   {\em expansion rate} and the {\em expansion volume rates}  of $\phi.$

 Let $K$ be a nonempty, compact and forward invariant set (i.e $\phi(K) \subset K$). The {\em expansion constant} of $\phi$ at $x$ is the positive number $$\tcb{E}(\phi,x) = \inf_{v \in T_x M \: \|v\|_x= 1} \|D\phi(x) v\|_{\phi(x)}.$$   Here $\|\cdot\|_x$ stands for the Riemaniann norm on $T_x M.$ Following Hirsch \cite{Hirsch94}, define the (logarithmic) {\em expansion rate} of $\phi$ at $K$ as
$${\cal E}(\phi, K) = \lim_{n \rar \infty} \frac{1}{n} \log(\min_{x \in K} E(\phi^n,x)),$$  where the limit exists by subadditivity.
The  {\em expansion rate} of $\phi$ is defined as $${\cal E}(\phi) = {\cal E}(\phi, M).$$
We let $\mathsf{Inv}(\phi)$ and $\mathsf{Inv}_{erg}(\phi)$ respectively denote the set of invariant (respectively ergodic) probability measures for $\phi.$

\tcb{Let $\mu \in \mathsf{Inv}_{erg}(\phi).$ By the Oseledec multiplicative ergodic theorem \cite{Os68}, there exist $k \in \{1, \ldots, d\},$ numbers $\Lambda^1 < \Lambda^2 < \ldots < \Lambda^k$   and, for $\mu$ almost of $x,$ vector spaces  $\{0\} = V_x^0 \subset V_x^1 \subset \ldots  V_x^{k}  = T_x M,$ such that for all $v \in V_x^{j} \setminus V_x^{j-1}$  $$\lim_{n \rar \infty} \frac{1}{n} \log \|D\phi^n(x) v\| = \Lambda^j.$$
The $\Lambda^j$ are called  the {\em Lyapunov exponents} of  $(\phi, \mu).$
The dimension of $(V_x^{j}), \, \mathsf{dim} (V_x^{j})$, depends only on $\mu$ and the number $\displaystyle \mathsf{dim} (V_x^{j}) - \mathsf{dim}(V_x^{j-1})$ is called the {\em multiplicity} of $\Lambda^j.$
We write $$\Lambda_1(\mu) \leq \ldots \leq \Lambda_d(\mu)$$ for the Lyapunov exponents of $(\phi, \mu)$ counted with their multiplicities.}

By a theorem of Schreiber \cite{Sch97},
\beq
\label{eq:Sch97}
{\cal E}(\phi) = \inf_{\mu \in \mathsf{Inv}_{erg}(\phi)} \Lambda_1(\mu). \eeq
%where $\Lambda_1(\mu)$ stands for the smallest Lyapunov exponent for $(\phi, \mu).$

For all $k \geq 0,$ we analogously define the {\em $k$-expansion volume rate} of $\Phi$ at $K$
as
$${\cal EV}_k(\phi,K) = \lim_{n \rar \infty} \frac{1}{n} (\min_{x \in K} \left [\log( J(\phi^n,x)) + k \log(E(\phi^n,x)) \right ]),$$ and the {\em $k$-expansion volume rate} of $\Phi$ as
\begin{equation}\label{eq:k-expansion volume rate}
{\cal EV}_k(\phi) = {\cal EV}_k(\phi,M).
\end{equation}
Again, these limits exist by subadditivity.

\tcr{Intuitively, the {\em expansion rate} measures the (asymptotic) rate   at which $\phi$ increases distance, and the {\em $0$-expansion volume rate } the (asymptotic) rate at which it increases volume. The $k$-expansion volume rate interpolates between these quantities.}

The following characterization easily follows from a beautiful result due to Schreiber \cite{Sch98} on the growth rates of sub-additive functions.
\bprop
\label{lem:EV}
The $k$-expansion volume rate of $\Phi$ is given by
\beq
\label{eq:EV}
{\cal EV}_k(\phi)  = \inf_{\mu \in \mathsf{Inv}_{erg}(\phi)} ((k+1) \Lambda_1(\mu) + \Lambda_2(\mu)
 + \ldots + \Lambda_d(\mu)), \eeq
where $\Lambda_1(\mu) \leq \ldots  \leq \Lambda_d(\mu)$ are the Lyapunov exponents of $(\phi,\mu)$ counted with their multiplicities.
\eprop
\prf
Let $F : M \times \NN \rar \RR$ be defined as $$F(x,n) = -\log J(\phi^n,x) -k E(\phi^n,x).$$ Then $F$ is continuous in $x$ and   subadditive with respect to $\phi,$ meaning that
$$F(x,n+1) \leq F(x,n) + F(\phi(x),1).$$ This directly follows from the properties  $J(\phi^{n+1},x) = J(\phi^n, \phi(x)) J(\phi,x)$ and $E(\phi^{n+1},x) \geq E(\phi^{n},\phi(x)) E(\phi,x).$
Therefore, by Theorem 1 in \cite{Sch98},
\begin{align*}
\lim_{n \rar \infty} \left(\sup_{x \in M} \frac{1}{n} F(x,n)\right) &= \inf_{n > 0}\left( \sup_{x \in M} \frac{1}{n} F(x,n)\right)\\  &= \sup_{\mu \in \mathsf{Inv}_{erg}(\phi)} \inf_{n > 0} \frac{1}{n} \int_M F(n,x) \mu(dx).
\end{align*}

For all  $\mu \in \mathsf{Inv}_{erg}(\phi)$ we have that
\begin{align*}
&\frac{1}{n} \int_M F(n,x) \mu(dx) \\
&= -\frac{1}{n} \sum_{k = 0}^{n-1} \int_M \log(J(\phi,\phi^k(x))) \mu(dx) - k  \frac{1}{n} \int_M \log(E(\phi^n,x)) \mu(dx)\\
&=- \int_M \log(J(\phi,x)) \mu(dx) - k  \frac{1}{n} \int_M \log(E(\phi^n,x)) \mu(dx).
\end{align*}
The first term on the right-hand side is equal to $-(\Lambda_1(\mu) + \ldots + \Lambda_d(\mu))$ by the multiplicative ergodic theorem \cite{Os68}, and the second term converges to $- k \Lambda_1(\mu).$
\qed
 \brem
 \label{rem:expans}
 {\rm We let $$\omega_{\phi}(x) = \bigcap_{n \geq 0} \clos{\{\phi^k(x)\:: k \geq n\}}$$ be the {\em omega limit set} of $x,$
 $$\displaystyle{\tcb{\mathsf{B}}(\phi) = \clos{\{x \in M \: : x \in \omega_{\phi}(x)\}}}$$ the {\em Birkhoff center} of $\phi$, and
 $$\displaystyle{\tcb{\mathsf{M}}(\phi) = \clos{\bigcup_{\mu \in \mathsf{Inv}_{erg}(\phi)} \mathsf{supp}(\mu)}}$$ the  {\em minimal center of attraction} of $\phi.$
 By the Poincaré recurrence theorem (see e.g.~\cite{Mane87}, Chapter 1), $\mathsf{M}(\phi) \subset \mathsf{B}(\phi).$
Thus, equalities (\ref{eq:Sch97}) and (\ref{eq:EV}) imply that $${\cal E}(\phi) = {\cal E}(\phi, \mathsf{B}(\phi))= {\cal E}(\phi, \mathsf{M}(\phi))$$ and
$${\cal EV}_k(\phi) = {\cal EV}_k(\phi, \mathsf{B}(\phi))= {\cal EV}_k(\phi, \mathsf{M}(\phi)).$$
These properties prove to be useful to compute or estimate the expansion and expansion volume rates in certain cases (see Examples \ref{ex:gradient} and \ref{ex:OnS2} below).
 }
 \erem
\brem
{\rm
\label{rem:boundsE}
We have that
$$d {\cal E}(\phi) \leq  {\cal EV}_0(\phi) \leq \log(\mathsf{deg}(\phi)).$$
The first inequality follows from identities (\ref{eq:Sch97}) and (\ref{eq:EV}), while the second follows from the second statement in the next proposition.

Note that this has the consequence that $$d {\cal E}(\phi) \leq  {\cal EV}_0(\phi) \leq 0$$ when $\phi$ is a diffeomorphism. Observe also that if ${\cal EV}_0(\phi) \leq 0,$ then $k \mapsto {\cal EV}_k(\phi)$ is nonincreasing.
}
\erem
We recall (see equation (\ref{eq:RLphi})) that for all $0 \leq k \leq r-1,$ ${\cal R}({\cal L}_{\phi}, C^{k}(M))$ is the spectral radius of ${\cal L}_{\phi}$ on $C^k(M).$
\bprop We have the following:
\label{prop:transfert}
\bdes
 \iti if  ${\cal  E}(\phi) > 0,$ then ${\cal R}({\cal L}_{\phi}, C^{r-1}(M)) = 1;$
 \itii if  ${\cal  E}(\phi) \leq  0,$ then  $$1 \leq {\cal R}({\cal L}_{\phi}, C^{r-1}(M)) \leq \mathsf{deg}(\phi) \max_{0 \leq k \leq r-1}e^{ -  {\cal  EV}_{k}(\phi)}.$$
 \edes
\eprop
\brem
\label{rem:campbell}
{\em
The first assertion of this proposition is a direct consequence of the seminal work of Ruelle (\cite{Ruelle89}). Some details are given below.

Some of Ruelle's results have been extended by Campbell and Latushkin in \cite{campbell} to the situation where $\phi$ is no longer expanding but is a covering map (i.e a local diffeomorphism as in the present setting). They compute the essential spectral radius of the transfer operator and provide an upper bound for the spectral radius in $C^0(M)$ (in the present setting) given by
\begin{equation}
\begin{split}\label{eq:Campbell-Latushkin}
&\exp{\Big( \sup_{\mu \in \mathsf{Inv}_{erg}(\phi)} \Big[H(\mu) - \int_M \log(J(\phi,x)) \mu(dx)\Big]\Big) }
\\= &\exp{\Big( - \inf_{\mu \in \mathsf{Inv}_{erg}(\phi)} \Big[(\Lambda_1(\mu) + \ldots +\Lambda_d(\mu)) - H(\mu)\Big]\Big)},
\end{split}
\end{equation}
where $H(\mu)$ is the measure-theoretic entropy of $(\phi,\mu).$
They claim (see \cite[Theorem 1]{campbell}) that this upper bound is also an upper bound for the spectral radius in $C^r(M)$ for $r \geq 1.$ Although this result is true when $\phi$ is expanding, it cannot be true when $\phi$ is not expanding, as shown by the following simple example. The error in their proof comes from the fact that they rely on  estimates (given in \cite{Ruelle89}) which are valid only for expanding maps.

The estimate given in Proposition \ref{prop:transfert}, (ii), provides a correct estimate  well-suited to non expanding maps.
}
\erem
\bex
\label{ex:countercampbell}
{\rm
We take $M = S^1 = \RR/\ZZ,$ and suppose that $\phi$ is a smooth, orientation preserving diffeomorphism with two fixed points, $0$ and $1/2$, such that $\phi$ coincides with  $$x \mapsto \frac{x}{\alpha}$$ on a neighbourhood of $0,$ \tcb{where} $ \alpha > 1$ and $\phi'(1/2) > 1.$  The ergodic measures of $\phi$ are the Dirac measures $\delta_0, \delta_{1/2}$, and for all $k  \geq 0,$ $${\cal EV}_k(\phi) =  - \ln(\alpha)(k+1) < 0.$$
Thus, by Proposition \ref{prop:transfert},  ${\cal R}({\cal L}_{\phi}, C^{r}(M)) \leq \alpha^{r+1}$ for all $0\leq r<\infty$. We now let $\rho(x) = \sin(2\pi x)$ if $r$ is odd, and $\rho(x) = \cos(2 \pi x)$ if $r$ is even. Then $$\|({\cal L}_{\phi^n}(\rho))\|_{C^{r}(M)} := \sum_{k = 0}^{r} \|({\cal L}_{\phi^n}(\rho))^{(k)}\|_0 \geq  |({\cal L}_{\phi^n}(\rho))^{(r)}(0)| =  \alpha^{ n (r+1)}.$$ This implies that  \begin{equation}\label{eq:spectral radius alpha r}{\cal R}({\cal L}_{\phi}, C^{r}(M)) =  \alpha^{r+1}\quad\text{for all}\quad 0\leq r<\infty.
\end{equation}
 This simple example shows that the inequality in Proposition \ref{prop:transfert} can be an equality, for any $r$.

The measure-theoretic entropy for any Dirac mass is $0$, whence we see that the Campbell-Latushkin upper bound in \eqref{eq:Campbell-Latushkin} is precisely $\alpha$. However, the authors claim in \cite[Theorem 1]{campbell} that this same upper bound for the $C^r$ spectral radius holds for all $0\leq r<\infty$, which cannot be true for any $r\geq 1$ by \eqref{eq:spectral radius alpha r}.
}
\eex
\subsubsection*{Proof of Proposition \ref{prop:transfert}}
{\em Step 1}. If ${\cal E}(\phi) > 0,$ then $\inf_{x \in M} E(\phi^n,x) \geq \theta > 1$ for some $n \geq 1$ and some $\theta > 1.$ Thus, replacing $\phi$ by $\phi^n$, we can assume that $EC(\phi,x) \geq \theta > 1.$ This condition means that $\phi$ is {\em expanding}. Then, by a theorem due to Ruelle \cite{Ruelle89}, Theorem 3.6 (ii) (see also \cite{Baladi2000}, Theorem 2.6), $R = {\cal R}({\cal L}_{\phi}, C^{r-1}(M))$ is an eigenvalue of ${\cal L}_{\phi}$ associated to a positive eigenfunction $\rho.$ Since $\int_M \rho dm = \int_M ({\cal L}_{\phi} \rho) dm,$ $R$ must be  $1.$ This proves the first assertion.
\medskip

\noindent {\em Step 2}. We now prove the left-hand side inequality of assertion $(ii).$ Suppose for the sake of contradiction that ${\cal R}({\cal L}_{\phi}, C^{r-1}(M)) < 1.$
Then $$\lim_{n \rar \infty} \|{\cal L}_{\phi}^n\|_{C^{r-1}(M)} = 0,$$ so that $\lim_{n \rar \infty} \|{\cal L}_{\phi}^n 1\|_{0} = 0$ in particular. On the other hand,
  $\int_M {\cal L}_{\phi}^n 1 dm =  \int_M 1dm = m(M) > 0.$ This is a contradiction.
\medskip

\noindent
{\em Step 3}. Our last  goal is to prove the  right-hand side inequality of assertion $(ii).$ It is convenient to firstly specify a norm on $C^{k}(M)$ for $k \geq 0.$

Throughout, $\RR^d$ is equipped with the Euclidean norm. For all $k \geq 1,$ let $L^k_{sym}(\RR^d)$ be the vector space of $k$-linear symmetric forms on $\RR^d.$ If $A : \RR^d \rar \RR^d$ is a linear map and $L \in L^k_{sym}(\RR^d), A^*L \in L^k_{sym}(\RR^d)$ is defined by $A^* L(u_1, \ldots, u_k) = L(Au_1, \ldots, Au_k).$
The norm of  $L \in L^k_{sym}(\RR^d)$ is defined as  $\|L\| = \sup \{|L(u_1, \ldots, u_k)| \: : u_i \in \RR^d, \|u_i\| \leq 1\}.$

We consider $U \subset \RR^d$ open and $f \in C^k(U) := \{ f : U \rar \RR, \, C^k\}.$ The $k$-th derivative  of $f$ is a continuous mapping $D^k f : U \rar L^k_{sym}(\RR^d).$ The following lemma will be used below. It follows \tcb{by induction} from classical rules in differential calculus.
%Its verification is left to the reader.
\blem
\label{lem:CDI}
Let $k \geq 1,$ and $U,V$ open subsets of $\RR^d.$
\bdes
\iti Let $g \in C^k(U).$ For all $f \in C^k(U)$ and $x \in U,$
$$\|D^k(g f)(x) - g(x) D^k f(x) \| \leq  \sum_{i = 0}^{k-1} \left(
                                                              \begin{array}{c}
                                                                k \\
                                                                i \\
                                                              \end{array}
                                                            \right)
 \|D^{k-i} g(x)\| \|D^i f(x)\|,$$ with the convention that $D^0 f = f.$
\itii Let $\Psi : U \rar V$ be a $C^k$ map. For all $f \in C^k(V)$ and $x \in U,$
\[
\|D^k(f \circ \Psi )(x) - D\Psi(x)^* D^k f(\Psi(x))\| \leq \]
\[  \sum_{i = 1}^{k-1} B_{k,i}(\|D\Psi(x)\|, \|D^2\Psi(x)\|, \ldots, \|D^{k-i+1}\Psi(x)\|)\|D^i f(\Psi(x))\|,
\]
where $(x_1, \ldots, x_{k-i+1}) \mapsto B_{k,i}(x_1, x_2, \ldots, x_{k-i+1})$ is a polynomial such that $B_{k,i}(x_1, 0, \ldots, 0)= 0.$
\edes
\elem
We now define a norm on $C^k(M).$ Let $W$ be the open ball in $\RR^d$ centered at the origin with radius $2$ and let $V$ be the open ball centered at the origin with radius $1.$

By the compactness of $M$ there exists an atlas $\{\alpha, {\cal O}_{\alpha}\}_{\alpha \in \aleph}$ with $\aleph$ finite such that:
\bdes
\iti $\alpha$ maps ${\cal O}_\alpha$ diffeomorphically onto  an open set in $\RR^d$ containing $\overline{W};$
\itii the open sets ${\cal O}'_{\alpha} = \alpha^{-1}(V), \alpha \in \aleph,$ cover $M.$
\edes
If $\rho  \in C^k(M)$  and  $1 \leq j \leq k,$ we set
 $$|\rho|_j = \sup_{\alpha \in \aleph, x \in \overline{V}}  \| D^j(\rho \circ \alpha^{-1})(x)\|$$ and
 \beq
 \label{eq:cknormatlas}
 \|\rho\|_{k} =  \|\rho\|_0 + \sum_{j = 1}^k |\rho|_j .
  \eeq
  It is not hard to verify that $\|.\|_{k}$ is a norm on $C^k(M)$ inducing the $C^k$ topology. For further reference we call this norm {\em the $C_k$ norm  induced by  $\{\alpha, {\cal O}_{\alpha}\}_{\alpha \in \aleph}.$}
\blem
\label{lem:spectralradius}
Let $k \geq 1$ and let $L : C^k(M) \rar C^k(M)$ be a bounded operator.
Suppose that there exist  sequences $(a_n)_{n\geq 0}, (b_n)_{n \geq 0}, a_n \geq 0, b_n \geq 0$  such that for all  $n \geq 0$ and  $\rho \in C^k(M)$,
$$|L^n \rho|_k  \leq  a_n |\rho|_k + b_n \|\rho\|_{k-1}.$$
Then $${\cal R}(L,C^k(M)) \leq \max \left({\cal R}(L,C^{k-1}(M)), \limsup_{n \rar \infty} a_n^{1/n} \right).$$
\elem
\prf
For all $\delta > 0,$ we set
$$ \|\rho\|_{k, \delta} = \|\rho\|_{k-1} + \delta |\rho|_k.$$ Note that
 $\|\rho\|_{k, \delta}$ and $\|\rho\|_{k}$ are equivalent norms. In particular we have
 $${\cal R}(L,C^{k}(M)) = \lim_{n \rar \infty} \|L^n\|_{k, \delta}^{1/n} \leq \|L\|_{k, \delta}$$
for all $\delta > 0.$

 We now fix $A >  \limsup_{n \rar \infty} a_n^{1/n}$ and
 $R > {\cal R}(L,C^{k-1}(M)).$  Then, for some $n \geq 0$ sufficiently large and all $\delta > 0,$
\begin{align*}
\|L^n \rho\|_{k, \delta} &\leq  \|L^n \rho\|_{k-1} + \delta [a_n |\rho|_k +  b_n \|\rho\|_{k-1}]\\
&\leq R^n \|\rho\|_{k-1} + \delta [A^n |\rho|_k +  b_n \|\rho\|_{k-1}]\\
&\leq \max \left(R^n + \delta b_n, A^n \right)  \|\rho\|_{k,\delta}.
\end{align*}
Thus $${\cal R}(L^n,C^{k}(M)) \leq \|L^n\|_{k, \delta} \leq \max \left(R^n + \delta b_n, A^n \right).$$
Since $\delta > 0$ is arbitrary, this shows that
$${\cal R}(L^n,C^{k}(M)) \leq \max \left(R^n, A^n \right).$$ Thus,
$${\cal R}(L,C^{k}(M)) = {\cal R}(L^n,C^{k}(M))^{1/n} \leq \max(A,R).$$ This concludes the proof. \qed

\blem We have the following:
\label{lem:spectralradius2}
\bdes
\iti
${\cal R}({\cal L}_{\phi}, C^{0}(M)) \leq \mathsf{deg}(\phi) e^{- {\cal  EV}_0(\phi)};$
\itii
for all $1 \leq k \leq r-1, {\cal L}_{\phi}$ satisfies the assumptions of Lemma \ref{lem:spectralradius} with
$$\limsup_{n \rar \infty} a_n^{1/n} \leq \mathsf{deg}(\phi) e^{ - {\cal  EV}_k(\phi)  }$$
\edes
\elem
\prf
Throughout this proof we set
\[
j_{\phi}( x) = \frac{1}{J(\phi,x)}.
\]
$(i):$ By the definition of ${\cal L}_{\phi},$  $$\|{\cal L}_{\phi} (\rho) \|_{0} \leq \mathsf{deg}(\phi) \sup_{x \in M} j_{\phi}(x)  \| \rho \|_{0}$$ for all $\rho \in C_0(M).$ Thus, replacing $\phi$ by $\phi^n$, we obtain that $$\|{\cal L}_{\phi}^n (\rho) \|_{0} \leq \mathsf{deg}(\phi)^n \sup_{x \in M} j_{\phi^n}(x)  \| \rho \|_{0},$$  whence the result follows from
the definition of ${\cal EV}_0(\phi).$

$(ii):$ To shorten notation we firstly consider the case where $\mathsf{deg}(\phi) = 1,$ so that $\phi$ is a diffeomorphism with inverse $\psi.$
Then ${\cal L}_{\phi} (\rho) =  (\rho \circ \psi)(j_{\phi} \circ \psi).$
Our first goal is to bound
$$|{\cal L}_{\phi} (\rho)|_k = \sup_{x \in \overline{V}, \alpha \in \aleph} \|D^k ({\cal L}_{\phi}(\rho)\circ \alpha^{-1})(x)\|.$$
We let $\alpha \in \aleph$ and $\overline{x} \in \overline{V}$, and choose $\beta \in \aleph$  such that $\psi(\alpha^{-1}(\overline{x})) \in {\cal O}'_{\beta}$ (recall that the family  $\{U'_{\beta}\}$ cover $M$).

 Set $U = \alpha (\psi^{-1}({\cal O}_\beta') \cap {\cal O}_\alpha),  f = \rho \circ \beta^{-1} : V \rar \RR, g = j_{\phi} \circ \psi \circ \alpha^{-1}  : U \rar \RR$ and
$\Psi = \beta \circ \psi \circ \alpha^{-1}: U \rar V.$
Then on $U$ we have $${\cal L}_{\phi}(\rho)\circ \alpha^{-1} = (f \circ \Psi) g.$$ Hence, relying on Lemma \ref{lem:CDI}, one can find a smaller neighbourhood of $\overline{x},$ $U_{\overline{x}} \subset U$, and a constant $C(\phi,\overline{x})$ (depending on $\phi$ and $\overline{x}$) such that for all $x \in U_{\overline{x}}$
\[
\begin{split}
&\|D^k ({\cal L}_{\phi}(\rho)\circ \alpha^{-1})(x) - g(x) D\Psi(x)^* D^k f(x)\|\\
&\leq C(\phi,\overline{x}) \left (|f(\Psi(x))| + \sum_{i = 1}^{k-1} \|Df^i(\Psi(x))\| \right ) \\
&\leq C(\phi,\overline{x}) \|\rho\|_{k-1}.
\end{split}
\]
We take constants $0<c,c'<\infty$ (depending only upon the atlas $\{\alpha, {\cal O}_\alpha\}$) such that for all $\alpha \in \aleph, x \in \alpha^{-1}(\overline{W})$ and $u \in T_x M$ we have $$c'\|u\|_x \leq \|D \alpha(x) u\| \leq c \|u\|_x.$$
Thus, defining $c'' = c/c',$ for all $x \in U_{\overline{x}}$ we have that
\[
\begin{split}
&\|g(x) D\Psi(x)^* D^k f(x)\| \leq g(x) \|D^k f(x)\| \|D\Psi(x)\|^k\\
&\leq c'' g(x) \|D^k f(x)\| \|D\psi(\alpha^{-1}(x))\|_{\alpha^{-1}(x)}^k \\
&=c''  j_{\phi}(\psi \circ \alpha^{-1}(x))\|D^k f(x)\| E(\phi,\psi \circ \alpha^{-1}(x))^{-k}\\
& \leq c'' \|D^k f(x)\| \sup_{y \in M} j(\phi,y) E(\phi,y)^{-k}.
\end{split}
\]
Finally, since $\overline{V}$ can be covered by finitely many neighbourhoods of the form $U_{\overline{x}},$
 we obtain that
 $$|{\cal L}_{\phi} (\rho)|_k \leq c'' |\rho|_k \sup_{y \in M} \Big[j(\phi,y) E(\phi,y)^{-k}\Big] + c_{\phi} \|f\|_{k-1},$$
 where $c''$ depends only upon the atlas $\{\alpha, {\cal O}_{\alpha}\}$ and $c_{\phi}$ depends on $\phi.$
 Replacing $\phi$ by $\phi^n$ gives
  $$|{\cal L}_{\phi^n} (\rho)|_k \leq c'' |\rho|_k \sup_{y \in M} \Big[j(\phi^n,y) E(\phi^n,y)^{-k}\Big] + c_{\phi^n} \|f\|_{k-1}.$$
This proves the desired result.

The proof for $\mathsf{deg}(\phi) > 1$ is similar, with the inverse of $\phi$ be replaced by the $\mathsf{deg}(\phi)$ local inverses.
 \qed
The proof of the right-hand side inequality of Proposition \ref{prop:transfert} $(ii)$ now easily follows from  lemmas \ref{lem:spectralradius} and \ref{lem:spectralradius2}.

\subsection{Application to random maps}
\label{sec:applRDS}
We recall that $P = P^{\nu}$, as defined in the beginning of the present section.
\bthm
\label{th:RDS3} We assume Assumption \ref{hyp:RDS1} and that $\nu_1 = \delta_{\phi}$ for some $\phi \in \mathsf{Diff_{loc}^r}(M).$
\bdes
\iti
If  ${\cal  E}(\phi) > 0,$ then   $\mathsf{Inv}(P) \subset {\cal M}_{ac}^{r-1}(M)$ for all $a < 1.$
\itii If  ${\cal  E}(\phi) \leq  0,$ then $\mathsf{Inv}(P) \subset {\cal M}_{ac}^{r-1}(M)$ for all
 $a  < \min_{k = 0, \ldots, r-1} \frac{e^{ {\cal  EV}_{k}(\phi)}}{\mathsf{deg}(\phi)}$.
 \edes
\ethm
\prf Theorem \ref{th:RDS3} follows from Theorem \ref{th:main2RDS} and Proposition \ref{prop:transfert} \qed
As an illustration of this last result, consider two examples where $\phi$ is a diffeomorphism, so that
$\min_{k = 0, \ldots, r-1} \frac{e^{ {\cal  EV}_{k}(\phi)}}{\mathsf{deg}(\phi)} = e^{{\cal  EV}_{r-1}(\phi)},$ and  where ${\cal  EV}_{r-1}(\phi)$ can be easily expressed.
\bex
\label{ex:gradient}
{\rm Suppose that $\phi$ is a $C^r$ diffeomorphism on $M$ such that for all $x \in M,$
 $$\omega_{\phi}(x) \subset \mathsf{Fix}(\phi):= \{p \in M \: : \phi(p) = p\}.$$ One can, for instance, imagine that $\phi = \Phi^1$ is the time one map of a flow $\{\Phi^t\}$ induced by  a $C^r, r \geq 1,$ gradient vector field $F = -\nabla V$ on $M$ (or more generally a vector field having a strict Lyapounov function).

Here $\mathsf{B}(\phi) = \mathsf{Fix}(\phi),$ so that by
Remark \ref{rem:expans}, $${\cal EV}_{r-1}(\phi) = {\cal EV}_{r-1}(\phi,\mathsf{Fix}(\phi)) = \inf_{p \in \mathsf{Fix}(\phi)} \log(J(\phi,p)) + (r-1) \Lambda_1(p)$$
 and $${\cal E}(\phi) = {\cal E}(\phi,\mathsf{Fix}(\phi)) = \inf_{p \in \mathsf{Fix}(F)} \Lambda_1(p).$$
 Here  $$J(\phi,p) = \log(|\mathsf{det} D\phi(p)|)$$ and $$\Lambda_1(p) = \min \{\log(|z|) \: : z \mbox{ is an eigenvalue of } D\phi(p)\}.$$

Note that, in case $\phi$ is the time one map of the flow induced by $F= -\nabla V,$ $\mathsf{Fix}(\phi) = \mathsf{Eq}(F) = F^{-1}(0)$, $J(\phi,p) = \mathsf{div}_p(F) = - \Delta V(p)$ and $\Lambda_1(p)$ is the smallest eigenvalue of the Hessian of $-V$ at $p.$

}
\eex

\bex
\label{ex:OnS2}
{\rm   We suppose here that $M = S^2$ and that $\phi = \Phi^1$ where $\{\Phi^t\}$ is induced by  a $C^r$ vector field $F.$ We no longer assume that $F$ is gradient-like but  will assume that $\mathsf{Eq}(F)$ is finite.

If $p \in  \mathsf{Eq}(F)$ we let  $$\Lambda_1(p) \leq \Lambda_2(p)$$ denote the real part
of the eigenvalues of $DF(p).$ Note that $$\mathsf{div}_p(F) = \Lambda_1(p) + \Lambda_2(p).$$
Given $T > 0,$ a $T$-{\em periodic orbit} is an orbit $\gamma = \{\Phi^t(p), t \in \RR\}$ such that $\Phi^T(p) = p$ and $\Phi^t(p) \neq p$ for all $0 < t < T$. We let $\mathsf{Per}_T(F)$ denote the set of such orbits and $\mathsf{Per}(F) = \cup_{T> 0} \mathsf{Per}_T(F).$

If $\gamma \in \mathsf{Per}_T(F)$ and $p \in \gamma,$  $D\Phi^T(p)$ has two (possibly equal) eigenvalues (that depend only on $\gamma$): $1$ (corresponding to the eigenvector $F(p)$) and  $J(\Phi^T,p).$   We let $$\{\Lambda_1(\gamma), \Lambda_2(\gamma)\} = \{0,\frac{\log(J(\Phi^T,p)}{T}\}$$ denote the logarithms of these eigenvalues, with the convention that $\Lambda_1(\gamma) \leq \Lambda_2(\gamma).$
A periodic orbit, $\gamma$, is said to be {\em linearly stable} if $\Lambda_1(\gamma) < 0.$ We let $\mathsf{Per}_{-}(F)$ denote the set of linearly stable periodic orbits. Note that, although  $\mathsf{Per}(F)$ may be uncountable, $\mathsf{Per}_{-}(F)$ is finite.

In the following lemma, Lemma \ref{lem:ergplan}, we implicitly identify an equilibrium point, $p$, with the orbit $\{p\} = \{\Phi^t(p) \: : t \in \RR\}.$
Again, combined with Theorem \ref{th:RDS3}, this gives simple conditions on $a$ ensuring the smoothness of invariant distributions.

\blem
\label{lem:ergplan}
Suppose that $F$ has finitely many equilibria.
Let $\mu$ be an ergodic probability measure for $\phi.$ Then
$\int \log (J(\phi,x)) \mu(dx) = \Lambda_1(\gamma) + \Lambda_2(\gamma)$ and $\Lambda_1(\mu) = \Lambda_1(\gamma)$ for some equilibrium or periodic orbit $\gamma.$ In particular, $${\cal EV}_{r-1}(\phi) = \min_{\gamma  \in \mathsf{Eq}(F) \cup \mathsf{Per}_{-}(F)} r \Lambda_1(\gamma) + \Lambda_2(\gamma)$$ and $${\cal E}(\phi) = \min_{\gamma  \in \mathsf{Eq}(F) \cup \mathsf{Per}_{-}(F)} \Lambda_1(\gamma).$$

\elem
\prf By the Poincaré recurrence theorem and the Birkhoff ergodic  theorem, there exists a set $\Omega \subset M,$ with $\mu(\Omega) = 1,$ such that $x \in \omega_{\phi}(x)$ (Poincaré) and $\frac{1}{n} \sum_{k = 0}^{n-1} \delta_{\phi^k(x)} \Rightarrow \mu$ (Birkhoff)  for all $x \in \Omega.$  Here $\Rightarrow$ stands for weak* convergence.

We take $p \in \Omega.$ We claim that $p$ is either a periodic point (i.e lies in a periodic orbit) or an equilibrium point for $\{\Phi^t\}.$
Clearly $\omega_{\phi}(p) \subset \omega_{\{\Phi^t\}}(p),$ the omega limit set of $p$ for $\{\Phi^t\}.$ Such a set is {\em internally chain recurrent} for $\{\Phi^t\}.$ Therefore, by a result proved in \cite{BH95}, Theorem 1.1, every point in $\omega_{\{\Phi^t\}}(p)$ is either periodic or belongs to an {\em orbit cycle}. An  orbit cycle is a finite sequence $\Gamma = {\gamma_1, \ldots, \gamma_m}$ of orbits such that the alpha limit set of $\gamma_i$ (for $\{\Phi^t\}$) is an equilibrium $e_{i-1}$ and its omega limit set is an equilibrium $e_i$, with $e_0 = e_m.$ Therefore, because $p \in \omega_{\{\Phi^t\}}(p),$  $p$ is either a periodic or an  equilibrium point. This proves the claim.

If $p$ is an equilibrium, then $\mu =  \delta_p,$  $\int \log(J(\phi, y)) \mu(dy) = \mathsf{div}_p(F) = \Lambda_1(p) + \Lambda_2(p)$, and $\Lambda_1(\mu) = \Lambda_1(p).$  If $p$ is $T$-periodic for $\{\Phi^t\}$ and $T = N/K$ is rational, then $\mu = \frac{1}{N} \sum_{i = 0}^{N-1} \delta_{\phi^i(p)}$ with $\phi^N(p) = \phi^{TK}(p) = p.$ Thus
\begin{align*}
&\int \log(J(\phi,y)) \mu(dy) = \frac{1}{TK} \log(J(\Phi^{T},p)^K) \\
&= \frac{1}{T}  \log(J(\Phi^T,p)) = \Lambda_1(\gamma) + \Lambda_2(\gamma).
\end{align*}
It $T$ is irrational, then $\mu = \frac{1}{T} \int_0^T \delta_{\Phi^s(p)} ds$ and again we have that
\begin{align*}
&\int \log(J(\phi,y)) \mu(dy)
 = \frac{1}{T}\int_0^T \log(J(\phi, \Phi^s(p))ds\\
& = \frac{1}{T} \int_0^T \int_0^1 \mathsf{Tr} (DF(\Phi^{s+u}(p))du ds= \frac{1}{T} \int_0^1 \int_0^T \mathsf{Tr} (DF(\Phi^{s+u}(p))du ds \\&=  \frac{1}{T} \int_0^T \mathsf{Tr} (DF(\Phi^{u}(p))du = \Lambda_1(\gamma) + \Lambda_{\tcb{2}}(\gamma).
\end{align*}
\qed
}
\eex
%%%%%%%%%%%%%%%
\section{Piecewise deterministic Markov processes}
\label{sec:pdmp}
We let $E$ be a finite set and $\{F_i\}_{i \in E}$ be a family of $C^r$ ($r \geq 1$) vector fields on $M$ where
$M$ is, as before, a $d$-dimensional compact connected Riemannian manifold.

We set $\M = M \times E.$ Then $\M$ can be viewed as a $d$-dimensional compact manifold with $\mathsf{card}(E)$ components. A map $g : \M \mapsto \RR$ is  $C^k$  if $x \rar g(x,i) = g_i(x)$  is $C^k$   for all $i \in E.$ A map $g: \M \mapsto \RR \cup \{\infty\}$ is lower semi-continuous if $g_i$ is lower semi-continuous for all $i\in E.$ The Riemannian measure on $\M$ is given by $\mathbf{m} = m \otimes \sum_{i \in E} \delta_i$, \tcb{where} $m$ is the Riemannian measure on $M$. The sets ${\cal M}_{ac}(\M), {\cal M}_{ac}^{ls}(\M)$ and ${\cal M}_{ac}^{r}(\M)$ are defined accordingly.

We let ${(Z_t = (X_t,I_t))_{t \geq 0}}$ be a continuous time Feller Markov process living on $\M$ whose  infinitesimal generator ${\cal A}$ acts on  functions $g \in C^1(\M)$ according to the formula
$${\cal A} g(x,i) = \langle F_i(x), \nabla g_i(x) \rangle_x + \sum_{j \in E} \alpha_{ij}(x)(g_j(x) -g_i(x)),$$
\tcb{where}:
\bdes
\iti  $\alpha_{ij}(x) \geq 0$ and (for convenience) $\alpha_{ii}(x) = 0$ for all $i,j\in E$;
\itii the matrix $(\alpha_{ij}(x))_{i,j \in E}$  is irreducible and  $C^{r-1}$ in $x.$
\edes
For further reference, we sometimes call the data $\{\{F_i\}_{i \in E}, (\alpha_{ij}(x))_{i,j \in E}\}$ the {\em characteristics} of $(Z_t)_{t \geq 0}.$

 An alternative pathwise description of the process is as follows. The component $(X_t)_{t \geq 0}$ is a solution to  the differential equation $$\frac{dX_t}{dt} = F_{I_t}(X_t),$$ while $(I_t)_{t \geq 0}$ is a jump process whose jump rates depends on $(X_t)$,
$$\Pr(I_{t+s} = j | \sigma(Z_u, u \leq t), I_t = i) = \alpha_{ij}(X_t) s + o(s).$$
In words, starting from $(x,i),$ $X_t$ follows the ODE induced by $F_i$ and switches to the ODE induced by $F_j$ at rate $\alpha_{ij}(X_t).$ Then $X_t$ follows the ODE induced by $F_j$ until it switches to the ODE induced by $F_k$ at rate $\alpha_{jk}(X_t)$, and so on.

This type of process falls under the broader category of {\em piecewise deterministic Markov processes}, introduced by Davis \cite{MR1283589}.
Their ergodic properties  have been the focus of much attention in the last decade (\cite{bakhtin&hurt}, \cite{BLMZ_2015} \cite{BLMZ_2019}, \cite{bakhtin&hurt&matt},  \cite{BCL17}, \cite{BHLM18}).
 \subsection{A discrete kernel associated to $(Z_t)_{t\geq 0}$}
In order to use the results of the preceding sections, we firstly introduce a (discrete time) Markov kernel $P$ whose invariant distributions are linked to the invariant distributions of $(Z_t)_{t\geq 0}.$

We let $\{\Phi^t_i\}_{t \in \RR}$ denote the flow induced by $F_i.$ We fix $\alpha > 0$ sufficiently large so that \tcr{for all $i \in E,$}
\beq
\label{eq:boundalpha0} \sup_{x \in M} \sum_{j \in E} \alpha_{ij}(x) < \alpha.
\eeq
Set $A_{ij}(x) = \frac{\alpha_{ij}(x)}{\alpha}$ for $i \neq j$ and
$A_{ii}(x) = 1 - \sum_{ j \neq i} A_{ij}(x).$  Let  $A, K$ and $P$ be the Markov operators on $\M$ respectively defined by
\beq
\label{eq:A}
A g(x,i) = \sum_j A_{ij}(x) g(x,j),
\eeq
\beq
\label{eq:K}
Kg(x,i) = \int_0^{\infty} \alpha e^{-\alpha t} g(\Phi_t^i(x),i) dt
\eeq and
\beq
\label{eq:KLambda}
P = K A
\eeq
\brem
\label{rem:chainXI}
{\rm The Kernel $P$ is the kernel of a discrete time chain $(X_n,I_n)_{n\geq 0}$  living on $\M$ whose dynamics can be described as follows. Starting from $(x,i) \in \M$, we pick a random variable $T$ having an exponential distribution with parameter $\alpha$, and set $X_1 = \Phi_i^T(x).$ We then choose $I_1 = j$ with probability $A_{ij}(X_1).$ }
\erem
Invariant distributions of the Markov kernel
$P$ and invariant distributions of the Markov process $(Z_t)_{t \geq 0}$ are linked by the following result proved in \cite[Proposition 2.4 and Lemma 2.6]{BLMZ_2015}.
\bprop
\label{prop:BLMZ24}
We let $(Z_t)_{t \geq 0}$ be the piecewise-deterministic Markov process having characteristics $\{\{F_i\}_{i \in E}, (\alpha_{ij}(x))_{i,j \in E}\}.$
The mapping $\mu \rar \mu K$ maps homeomorphically  $\mathsf{Inv}(P)$ (respectively  $\mathsf{Inv}_{erg}(P),$ the set of ergodic probability measures of $P$) onto the set of invariant (respectively ergodic) probability measures for $(Z_t)_{t \geq 0}.$ Its \tcb{inverse} homeomorphism is given by $\mu \mapsto \mu A.$

Moreover we have that $\mathsf{supp}(\mu) = \mathsf{supp}(\mu K)$ for all $\mu \in \mathsf{Inv}(P).$
\eprop
 By Liouville's formula, the transfer operator of $\Phi^t_i$ (see Section \ref{sec:RDS}) is given by
\beq
\label{eq:LPhit}
{\cal L}_{\Phi_i^t}(\rho)(x) = \rho(\Phi_i^{-t}(x)) \exp{[-\int_0^t \mathsf{div}(F_i)(\Phi^{-s}_i(x)) ds]}
 \eeq for  $\rho \in L^1(m),$ where $\mathsf{div}(F_i)$ denotes the divergence of $F_i$ on $M.$
 We also set
 \beq
\label{eq:Li}
{\cal L}_{i}(\rho)(x) = \int_0^{\infty} \alpha e^{-\alpha t} {\cal L}_{\Phi_i^t}(\rho)(x) dt
 \eeq
for  $\rho \in L^1(m)$. \tcr{This integral is well defined, as the integral of a nonnegative function, but may be infinite for small values of $\alpha.$ However, it  is always finite for $\alpha$ sufficiently large (see Lemma \ref{lem:KLambda}, (iii)).} Observe that, using  the notation of Proposition \ref{lem:Lnu1},  ${\cal L}_i :={\cal L}_{\nu_1}$ where $\nu_1$ is the measure on $\mathsf{diff}^r_{loc}(M)$ given by
$\nu_1 = \int_0^{\infty} \alpha e^{-\alpha t} \delta_{\Phi^t_i} dt.$

Associated to $K$ is the transfer operator defined on $L^1(\m)$ by
$${\cal K} \rho(x,i) = {\cal L}_i \rho_i(x).$$
The purpose of the next lemma is twofold. Firstly, it will be used to show that $P$ satisfies Assumption \ref{hyp1}, $(iii),$  with $\C(\M)$ one of the sets ${\cal M}_{ac}(\M), {\cal M}_{ac}^{ls}(\M)$ or ${\cal M}_{ac}^r(\M).$ Secondly, it shows that the mapping $\mu \rar \mu K$ in Proposition \ref{prop:BLMZ24} preserves these sets.
\blem
\label{lem:KLambda}
Suppose that $\mu \in {\cal M}_{ac}(\M)$ has  density $\rho$ with respect to $\mathbf{m}.$
Then we have the following:
\bdes
\iti  $\mu A$ has density $ A^t \rho$ given by
$$A^t \rho (x,i) = \sum_j \rho_j(x)A_{ji}(x).$$ If $\rho$ is lower semi-continuous or $C^k$ with $0 \leq k \leq r-1,$ then so is $A^t \rho.$
\itii $\mu K$ has a density given by ${\cal K} \rho.$ If $\rho$ is lower semi-continuous, then so is ${\cal K} \rho.$
\itiii If we furthermore assume that
\beq
\label{eq:boundalpha1}
\alpha > \max_{i \in E} \log \left({\cal R}({\cal L}_{\Phi_i^1}, C^{r-1}(M))\right),
  \eeq
  then ${\cal K}$ is a bounded operator on $C^{r-1}(\M)$ and
$${\cal R}({\cal K}, C^{r-1}(\M)) \leq \frac{\alpha}{\alpha -  \max_{i \in E} \log \left({\cal R}({\cal L}_{\Phi_i^1}, C^{r-1}(M))\right)}.$$
\edes
\elem
\prf $(i)$ is immediate to verify and $(ii)$ easily follows from  Proposition \ref{lem:Lnu1}.

We now turn to $(iii)$. By classical results (see \cite[Chapter V, Corollary 4.1]{Hartman82} for example), $(t,x) \mapsto \Phi^t_i(x)$ is $C^r.$
The form of ${\cal L}_{\Phi^t_i}$ (see equation (\ref{eq:LPhit})) and  the fact that $\mathsf{div}(F_i)$ is $C^{r-1}$ imply that
$$\sup_{0 \leq t \leq 1} \|{\cal L}_{\Phi^t_i}\|_{C^{r-1}(M)} \leq  C$$ for some constant $C <\infty$ (depending on $r$).
For $t \geq 0,$ we write $t = n + s$ for $n \in \NN$ and  $0 \leq s \leq 1.$ Thus
$${\cal L}_{\Phi^t_i} = {\cal L}_{\Phi^1_i}^n \circ {\cal L}_{\Phi^s_i}.$$ Therefore for all $\eps > 0$ there exists another constant $C'<\infty$ such that for all $t \geq 0$ we have
\beq
\label{eq:boundnormL}
\|{\cal L}_{\Phi^t_i}\|_{C^{r-1}(M)} \leq C' e^{ n(\log(R_i) + \eps)} \leq C' e^{t(\log(R_i) + \eps)},
\eeq where $R_i$ stands for
${\cal R}({\cal L}_{\Phi_i^1}, C^{r-1}(M)).$
  Proposition \ref{lem:Lnu1} then implies that ${\cal L}_i$ is a bounded operator on $C^{r-1}(M).$ We likewise have that ${\cal K}$ is a bounded operator on $C^{r-1}(\M).$

  We now establish the upper bound on the spectral radius. Note that for all $n \in \NN$ we have
$${\cal K}^n \rho(x,i) = \mathbb{E}({\cal L}_{\Phi_i^{S_n}}(\rho_i)(x)),$$ where $S_n = T_1 + \ldots + T_n$ and $\{T_i\}_{i \geq 1}$ is a sequence of independent random variables having an exponential distribution with parameter $\alpha.$
Thus
$$\|{\cal K}^n \rho\|_{C^{r-1}(M)} \leq \max_{i \in E} C' \mathbb{E}[e^{S_n(\log(R_i) + \eps)}] =  \max_{i \in E} C' \big(\mathbb{E}[e^{T_1(\log(R_i) + \eps)}]\big)^n.$$
This proves the result.
\qed
\subsection{Invariant distributions}
Let $C_{pc}(\Rp,E)$ be the set of piecewise continuous functions $J : \Rp \rar E.$ Given $J \in C_{pc}(\Rp,E),$ we let $t \rar \Phi^t(x,J)$ denote the solution to the non-autonomous differential equation
\beq
\label{eq:control}
\frac{dx}{dt} = F_{J(t)}(x), \eeq
with initial condition $x(0) = x.$ For all $x \in M,$ we define
$$\gamma^+(x) =  \{\Phi^t(x,J) \: : t \geq 0 \mbox{ and } J \in C_{pc}(\Rp,E) \}.$$ We let $\Gamma$ be the possibly empty, compact connected set defined by $$\Gamma = \bigcap_{x \in M} \clos{\gamma^+(x)}.$$
Connectedness (as well as other topological properties of $\Gamma$) are proved in \cite[Proposition 3.11]{BLMZ_2015} (see also the erratum \cite{BLMZ_2019}).
By Proposition 3.13 in \cite{BLMZ_2015} we have $$\Gamma_P = \Gamma \times E,$$ where $\Gamma_P$ is  {\em the accessible set} (as defined in Section \ref{sec:onH1H2}) of the kernel $P$ given by (\ref{eq:KLambda}).

We let $r_{max} \in \{1, 2, \ldots \} \cup \{\infty\}$ be the maximal $r$ such that all the $F_is$ are $C^{r}.$ We define $\textbf{F}_0:= \{F_i \:  : i \in E\}$ and inductively, for all $n = 1, \ldots r_{\max}-1,$ $\textbf{F}_{n} = \textbf{F}_{n-1} \cup \{[F,G] : F  \in \textbf{F}_0, G \in \textbf{F}_{n-1}\}$, where $[F,G]$ is the Lie bracket of $F$ and $G.$

We let $n \leq r_{max}-1.$ Inspired by  the terminology used in \cite{BLMZ_2015} (see also \cite[Chapter 6]{BenaimHurth}), we say that a point $p \in M$ satisfies the $n$-{\em weak bracket condition} if  $\textbf{F}_{n}(p):= \{G(p) \: : G \in \textbf{F}_n\}$ spans $T_p M.$
We say that $p$ satisfies the {\em weak bracket condition} if it satisfies the $n$-weak bracket condition, for some $n \leq r_{max}-1.$

It was proved  in \cite{bakhtin&hurt} (for $\alpha_{ij}(x)$ constant over $x$) and in \cite{BLMZ_2015} that for $C^{\infty}$ vector fields (i.e $r_{max} = \infty$), the existence of a point $p \in \Gamma$ at which the weak bracket condition holds implies that $(Z_t)$ has a unique invariant distribution which is absolutely continuous with respect to $\m.$ The next theorem also shows that its density is lower semi-continuous. A first version of this result, when $\alpha_{ij}(x)$ is constant over $x$, was proved in \cite{BT23}.

\bthm
\label{th:lscpdmp0} Assume there exists a point $p \in \Gamma$ at which the weak bracket condition holds. Then $(Z_t)$ has a unique invariant probability measure $\Pi$ which is absolutely continuous with respect to $\m$ and whose density $\rho$ is lower semi-continuous. In addition:
\bdes
\iti $\mathsf{supp}(\Pi) = \Gamma \times E;$
\itii For all $i \in E,$
 \beq
 \label{eq:supprhoi}
 \mathsf{supp}(\rho_i):= \clos{\{x \in M: \: \rho_i(x) > 0\}} = \Gamma, \eeq
and
\beq
\label{eq:intersupprhoi}
\bigcup_{t \geq 0} \Phi^t_i\left(\mathsf{Int}(\Gamma) \cap {\cal WB}(M)\right) \subset \{x \in M: \: \rho_i(x) > 0\} \subset \mathsf{Int}(\Gamma),\eeq where
${\cal WB}(M)$ stands for the open set of points at which the weak bracket condition holds.
 \edes
\ethm
\prf
We let $(p,i_0) \in \Gamma \times E = \Gamma_P.$ By Theorems 4.1 and 4.4 in \cite{BLMZ_2015}, $(p,i_0)$ is  a weak Doeblin point (as defined in Section \ref{sec:onH1H2}) of $P$ with a minorizing measure given by
\beq
\label{eq:doeblinatp}
\pi(dx di)
= c \Ind_{{\cal V} \times E}(x,i) \m(dx di),
  \eeq for some   nonempty open set ${\cal V} \subset M$   and $c > 0.$

Thus, $\pi(dx di)  \geq \theta(x) \m(dx di) : = \pi'(dx di)$ where $0 \leq \theta \leq c \Ind_{{\cal V}}$ is  continuous and nonzero somewhere. In particular,$\pi' \in {\cal M}_{ac}^{ls}(\M)$  and $\mu - \pi' \in {\cal M}_{ac}^{ls}(\M)$ for all $\mu \in {\cal M}_{ac}^{ls}(\M)$ greater than $\pi'.$ Therefore, by Theorem \ref{th:piirreducible}, $P$ has a unique invariant distribution $\mu$ having a lower semi-continuous density $h.$ By Proposition \ref{prop:BLMZ24} and Lemma \ref{lem:KLambda}, $\Pi = \mu K$ is the unique invariant distribution of $(Z_t)_{t\geq 0}$ and its density, $\rho = {\cal K} h,$ is lower semi-continuous. Also $\mu$ and $\Pi$ have the same support.

Basic properties of the accessible set (see \cite[Proposition 5.8 (iv)]{BenaimHurth}, for example) imply that $\mathsf{supp}(\mu)$ (hence $\mathsf{supp}(\Pi)$) is equal to $\Gamma_P.$
Clearly $\mathsf{supp}(\Pi) \subset \mathsf{supp}(\rho).$
Conversely, if $\rho_i(x) > \delta > 0,$ by the lower semi-continuity of $\rho,$ there exists a ball $B(x,\eps)$ such that $\rho_i(y) > \delta$ for all $y \in B(x,\eps).$
Thus $\Pi(B(x,\eps) \times \{i\}) > 0.$  This proves the converse inclusion $\mathsf{supp}(\rho) \subset \mathsf{supp}(\Pi).$

%\tcr
{
We now turn to the proof of \ref{eq:intersupprhoi}. The right hand side inclusion is immediate because $\{\rho_i > 0\}$ is an open set by lower semi-continuity and is contained in $\Gamma$ by what precedes.

For the left hand side inclusion, let $q \in \mathsf{Int}(\Gamma) \cap {\cal WB}(M).$ We claim that there exists $\delta > 0$ and a neighborhood of $q,$ ${\cal V} \subset \mathsf{Int}(\Gamma) \cap {\cal WB}(M)$ such that the density $h$ of $\mu$ (the invariant probability measure of $P$) is almost surely bounded below by $\delta$ on ${\cal V} \times E.$ The proof of the this claim is inspired by the proof of Lemma 8.1 in \cite{BT23} but we give details for the sake of clarity.

Because, the weak Bracket condition holds at $q$ for the family $\{F_i\}_{i\in E},$
it also holds for the family $\{-F_i\}_{i\in E}.$  Therefore, by Lemma 6.19 in \cite{BenaimHurth} (see also the proof of Theorem 4.4 in \cite{BLMZ_2015}), for every $\eps >0$ there exists $\i = (i_1, \ldots, i_d) \in E^d$ and
$t^* = (t_1^*, \ldots, t_d^*) \in (0,\eps)^d$ such that the mapping
$$t \in (0,\eps)^d \mapsto \mathbf{\Psi}^t(q):= \Phi_{i_1}^{-t_1} \circ \ldots \circ \Phi_{i_d}^{-t_d}(q)$$ is a submersion at $t^*.$
Let $$\mathbf{\Phi}^t = (\mathbf{\Psi}^t)^{-1} = \Phi_{i_d}^{t_d} \circ \ldots \circ \Phi_{i_1}^{t_1}.$$
Then, $\mathbf{\Phi}^t \circ \mathbf{\Psi}^t(q) = q$ and, consequently,
$$\partial_t \mathbf{\Phi}^{t^*}(p) + D \mathbf{\Phi}^{t^*}(p) \partial_t \mathbf{\Psi}^{t^*}(q) = 0,$$
where $p = \mathbf{\Psi}^{t^*}(q).$ This shows that
$t \mapsto \mathbf{\Phi}^t(p)$ is a submersion at $t^*.$ because $D \mathbf{\Phi}^{t^*}(p)$ is invertible and $\partial_t\mathbf{\Psi}^{t^*}(q)$
has rank $d.$ Also, by choosing $\eps$ sufficiently small, we can always assume that $p \in \mathsf{Int}(\Gamma).$ By Theorem 4.1 in \cite{BLMZ_2015}, this
implies that $(p,i_0)$ (for any $i_0 \in E$) is a weak accessible Doeblin point of $P$ with a minorizing measure $\pi$ given by (\ref{eq:doeblinatp}) where
${\cal V}$ is a {\bf neighborhood of $q.$} Therefore, (proceeding like in the proof of Theorem \ref{th:piirreducible}), one has
 $R_a^2((x,j), \cdot) \geq \delta \pi(\cdot)$ for some $\delta > 0$ and all $(x,j).$ Thus $\mu \geq \delta \pi$ with $\mu$ the invariant probability measure of $P.$ The density $h$ of $\mu$ is therefore almost surely bounded below by $\delta$ on  ${\cal V} \times E.$ This proves the claim.

 Let now $(U_n)_{n \geq 1}$ be a family of open set such that  $U_n \subset \overline{U_n} \subset U_{n+1} \subset  \mathsf{Int}(\Gamma) \cap {\cal WB}(M)$ and $\cup_{n \geq 1} U_n = \mathsf{Int}(\Gamma) \cap {\cal WB}(M).$ Covering $\overline{U_n}$ by a finite family of open sets $\{{\cal V}\}$ that are like in the claim, we get that $h$ is almost surely bounded below by some $\delta_n > 0$ on $U_n.$ Set
 $h_n = \max(h, \delta_1 \Ind_{U_1}, \ldots, \delta_n \Ind_{U_n}).$ Then $h_n$ is l.s.c, $h_n = h$ almost surely on $U_n$ and $h_n = h$ on $M \setminus U_n.$ Also, $h_n$ converges, as $n \rar \infty,$ toward a $l.s.c$ function $\tilde{h}$ which is positive on $\mathsf{Int}(\Gamma) \cap {\cal WB}(M)$ and which equals $h$ almost surely. The inclusion to be proved then follows from the identity $\rho_i(x) = {\cal K}(h)(x,i).$}
 \qed
%%%%%%%%%%%%%
The next result considers the situation where $\Gamma$ is empty but the weak bracket condition holds everywhere. It relies on the preceding result combined with ideas and results from \cite{BCL17}.
\bthm
\label{th:lscpdmp1}
We assume that the weak bracket condition holds at every point $p \in M.$ Then $(Z_t)_{t\geq 0}$ has finitely many   ergodic probability measures $\Pi^1, \ldots, \Pi^k$. These are absolutely continuous with respect to $\m$, with lower semi-continuous densities $\rho^1, \ldots, \rho^k$. For each $j = 1, \ldots, k$, the support of $\Pi^j$ can be written as $\mathsf{supp}(\Pi^j) = \Gamma^j \times E$, where $\Gamma^j$ is a compact connected set. Furthermore, for all $i \in E,$ $$\mathsf{supp}(\rho_i^j):= \clos{\{x \in M: \: \rho_i^j(x) > 0\}} = \Gamma^j$$
and $$\{x \in M: \: \rho_i^j(x) > 0\} = \mathsf{Int}(\Gamma^j).$$
\ethm
\prf The proof uses some results and ideas from control theory. For consistency with the terminology used in \cite{BLMZ_2015}, we phrase it using differential inclusions. We let $$\mathsf{co}(F)(x) = \Big\{\sum_{i \in E} p_i F_i(x)\: : p_i \geq 0, \sum_{i \in E} p_i = 1\Big\} \in T_x M$$ be the convex hull of the family $\{F_i(x)\}_{i \in E}.$
A solution to the differential inclusion
\beq
\label{eq:DI}
\dot{\eta} \in \mathsf{co}(F)(\eta)
 \eeq is an absolutely continuous function $\eta \in C^0(\RR_{\geq 0}, M)$ which satisfies $\dot{\eta}(t) \in \mathsf{co}(F)(\eta(t))$ for almost all $t \in \RR_{\geq 0}.$ Such a differential  inclusion induces a set-valued dynamical system defined as $$\Psi_t(x) = \{\eta(t) : \eta(0)= x \mbox{ and } \eta \mbox{ is solution to } (\ref{eq:DI})\}.$$ We refer the reader to \cite{BLMZ_2015} for background and references.
 For $I \subset \RR,$ we set $\Psi_I(x) = \bigcup_{t \in I} \Psi_t(x).$
 We call a set $C \subset M,$ a {\em compact invariant control set} if $C$ is nonempty, compact and $C = \overline{\Psi_{[0,\infty)}(x)}$ for all $x \in C.$ This is consistent with the terminology used in control theory (see, for instance, \cite[Definition 2.4  and Theorem 2.2]{BCL17}).
 The set $\Gamma$ previously defined is, when it exists, a compact invariant control set. This follows, for instance, from \cite[Proposition 3.11]{BLMZ_2015}.
 Under the present assumption that the weak bracket conditions holds at every point $p \in M$, there are, by \cite[Corollary 2.13]{BCL17}, finitely many compact invariant control sets $\Gamma^1, \ldots, \Gamma^k.$ Furthermore we have the following:
 \bdes
 \iti for all $j \in \{1, \ldots, k\}$ $\clos{\mathsf{Int}(\Gamma^j)} = \Gamma^j;$
 \itii for each $x \in M,$ there exists $j \in \{1, \ldots, k\}$ such that $\gamma^+(x) \cap \mathsf{Int}(\Gamma^j) \neq \emptyset;$
 \itiii for each $j \in \{1, \ldots, k\}$ and $x \in \Gamma^j, \mathsf{Int}(\Gamma^j) \subset \gamma^+(x).$
 \edes
 It follows from $(i), (iii)$ and the definition of a compact invariant control set, that
  $\Gamma^j = \bigcap_{x \in \Gamma^j} \clos{\gamma^+(x)}.$ The proof of  Theorem \ref{th:lscpdmp1} then applies verbatim to $P$ restricted to $\Gamma^j.$ This proves that $P$ restricted to $\Gamma^j$ has a unique, hence ergodic for $P$, invariant distribution $\Pi^j$ with density $\rho^j$ enjoying the properties stated in the theorem.

To establish that the $\Pi^j$s are the only ergodic probability measures, it suffices to show that every $\mu \in \mathsf{Inv}(P)$ is supported on $\bigcup_{j = 1}^k \Gamma^j.$ It easily follows from $(ii)$ that  $W = \bigcup_{j = 1}^k \mathsf{Int}(\Gamma^j)$ is accessible for $P$, that is $R_a(x,W) > 0$ for all $x \in M$ (this can, for instance, be deduced from the support theorem, \cite[Theorem 3.4]{BLMZ_2015,BLMZ_2019}). By the Feller continuity of $R_a$ (inherited from the Feller continuity of $P$), the Portmanteau theorem and the compactness of $M,$ we have that $R_a(x,W) \geq \delta > 0$ for all $x \in M$, for some $\delta > 0.$
Since $R_a(y,M \setminus \clos{W}) = 0$ for all $y \in \clos{W}$ one obtains that (one may compare this to \cite[Theorem 4.7]{BCL17})
\begin{align*}
&\mu(M \setminus \clos{W}) = \mu R_a^2(M \setminus \clos{W}) = \int_{M \setminus \clos{W}} \mu R_a(x,dy) R_a(y, M\setminus \clos{W})\\
&\leq (1-\delta) \mu R_a(M \setminus \clos{W}) = (1-\delta) \mu(M \setminus \clos{W}).
\end{align*}
We therefore obtain that $\mu(M \setminus \clos{W}) = 0.$\qed
\subsection{Smooth invariant distributions on the torus}
\label{sec:torus}
%%%%%%%%%%%%%%%%%%%%%%%%%%%%%%%%%
This section is motivated by the work of Bakhtin, Hurth, Lawley and Mattingly \cite{BHLM18}. It retrieves and substantially extends their main result (see Remark \ref{rem:torus}).

Here we assume that $M = \mathbb{T}^2 = \RR^2/\ZZ^2$ is the two dimensional flat torus, $E = \{1,2\}$, and that the vector fields $F_1, F_2$ are $C^r$ with  $r \geq 2,$ and transverse everywhere - that is $\{F_1(p), F_2(p)\}$ span $T_p \mathbb{T}^2$ for all $p.$
In particular $F_1, F_2$ never vanish.
Moreover we assume that the jump rates are constant, that is $$\alpha_{12}(x) = \alpha_{12} > 0, \mbox{ and } \alpha_{21}(x) = \alpha_{21} > 0.$$

Using the notation introduced in Example \ref{ex:OnS2}, we let  $\mathsf{Per}_{-}(F_i)$ denote the (possibly empty) finite set of  linearly stable periodic orbits of $F_i.$
For $\gamma \in \mathsf{Per}_{-}(F_i)$ we let $\Lambda_{1,i}(\gamma) < 0$ denote the non-zero Floquet exponent of $\gamma.$

We shall establish here the following result.

 \bthm
\label{th:torus} We let $1 \leq k \leq r.$
Assume  that for all  $i = 1, 2$ and $\gamma \in \mathsf{Per}_{-}(F_i)$,
$$\min (\alpha_{12}, \alpha_{21}) > - k \Lambda_{1,i}(\gamma),$$
\tcr{with the convention that the left-hand side is  zero when $\mathsf{Per}_{-}(F_i) = \emptyset.$}
Then $(Z_t)$ has finitely many ergodic probability measures (see Theorem \ref{th:lscpdmp1}), each of which has  a $C^{k-1}$ density with respect to $\m.$
\ethm
\bcor
\label{cor:torus}
Suppose that $F_1$ \tcr{has no periodic orbit} and that $F_2$ has no \tcr{linearly stable} periodic orbit. Then $(Z_t)$ has a unique invariant distribution and its density is $C^{r-1}.$
\ecor
\prf A fixed-point-free $C^2$ flow with no periodic orbits on  $\mathbb{T}^2$ has dense orbits (see the proof of Proposition \ref{prop:exptorus}). The accessible set is then $\mathbb{T}^2$ and uniqueness follows (see e.g~Theorem \ref{th:lscpdmp0}). The $C^{r-1}$ continuity follows from Theorem \ref{th:torus}. \qed
\brem
\label{rem:torus}
{\rm
Using ideas inspired  by Malliavin calculus, Bakhtin, Hurth, Lawley and Mattingly gave in \cite{BHLM18} a proof of Corollary \ref{cor:torus}, (when  $r = \infty$ and  $\alpha_{12} = \alpha_{21}$) in the particular case where each of the flows induced by $F^1$ and $F^2$ possess an invariant probability measure with an everywhere  positive $C^{\infty}$ density. This, it should be noted, is a strong assumption.}
\erem

\subsubsection*{Proof of Theorem \ref{th:torus}}
The idea of the proof is to show that $P^n$ (for $n$ sufficiently large) satisfies the standing assumption, Assumption \ref{hyp1}, and the assumptions of Theorem \ref{th:main2}.
We assume here that $F_1, F_2$ are $C^r$ with $r \geq 1.$ The assumption that $r \geq 2$ will be required in Proposition \ref{prop:exptorus}.

We let $(X_n,I_n)_{n \geq 0}$ be the discrete-time Markov chain with kernel $P$ (see Remark \ref{rem:chainXI}), and define $\tau = \min\{k \geq 1: \: I_k \neq I_0\}$ to be the first switching time.
For $n \geq 2$ and $1 \leq k \leq n-1$
we set $$P_{n,k}(f)(x,i) =  \mathbb{E} \left[f(X_n,I_n) \Ind_{\tau = k} | (X_0, I_0 ) = (x,i)\right ]$$ and
$$\Delta_{n,n}(f)(x,i) =  \mathbb{E} \left[f(X_n,I_n) \Ind_{\tau \geq n} | (X_0, I_0 ) = (x,i) \right ].$$
Clearly we have that
$$P^n f = \sum_{k = 1}^{n-1} P_{n,k} f + \Delta_{n,n} f.$$
We now decompose the (matrix) operator $A$  as
$A = S + \bar{S}$ where $S$ (corresponding to switching) is defined by
$$Sf(x,i) = A_{ij} f(x,j) \mbox{ with } j = 3-i,$$ and $\bar{S}$ (corresponding to not switching) is given by
$\bar{S} f(x,i) = A_{ii} f(x,i).$
It is readily seen that
$$P_{n,k}  = (K \bar{S})^{k-1} KS P^{n-k}   = (K \bar{S})^{k-1} [KSK] A P^{n-k-1}.$$
This simply express the fact that the first switch occurs at time $k.$
We likewise have that
$$\Delta_{n,n} = (K \bar{S})^{n-1} P.$$

In the next three lemmas, we use the following convenient notation. We denote $\C(\M) = {\cal M}_{ac}^{r-1}(\M),$ and if $\mu \in \C(\M)$ has density $\rho,$ then $\|\rho\|_{C^{r-1}(\M)}$ is denoted by $\|\mu\|_{\C(\M)}$. We also assume that the parameter $\alpha$ that occurs in the definitions of $A$ and $K$ satisfies inequality (\ref{eq:boundalpha1}). \tcr{That is $$\alpha> \max_{i = 1,2} \log \left({\cal R}({\cal L}_{\Phi_i^1}, C^{r-1}(M))\right).$$}

\tcr{The next lemma simply expresses the fact that "switching creates density".}
\blem
\label{lem:KSK} We suppose that $F_1, F_2$ are transverse at  every point $p \in \mathbb{T}^2.$ For all $\eps > 0,$ $KSK$ can be decomposed as $KSK = Q + \Delta$ where $Q, \Delta$ are Feller sub-Markov kernels and satisfy:
\bdes
\iti ${\cal M}(\M) Q \subset \C(\M)$;
\itii $\C(\M) \Delta \subset \C(\M)$;
\itiii $\|\mu \Delta\|_{\C(\M)} \leq \eps \|\mu\|_{\C(\M)}$ for all $\mu \in \C(\M)$.
\edes
\elem
\prf We set $j = 3-i$ for $i \in \{1,2\}.$ We note that
$$KSK f(x,i) = A_{ij} \int_{\Rp^2} f(\Phi_j^t \circ \Phi^s_i(x),j) \alpha^2 e^{-\alpha (t+s)} dt ds.$$
For all $n > 1,$ we let $\eta_n : \RR \mapsto \Rp$
be a $C^{\infty}$ function such that $\eta_n = 1$ on $[\frac{1}{n}, n],$ $\eta_n = 0$ on $\RR \setminus [\frac{1}{2n}, 2n]$, and $0 \leq \eta_n\leq 1.$ We set
 $$Q f(x,i) = A_{ij} \int_{\Rp^2} f(\Phi_j^t \circ \Phi^s_i(x),j) \alpha^2 e^{-\alpha (t+s)} \eta_n(t) \eta_n(s) dt ds$$ and $\Delta = KSK - Q.$
The assumption that $F_1, F_2$ are transverse makes the map $(t,s) \in \Rp^2 \rar \Phi_j^t \circ \Phi^s_i(x) \in \mathbb{T}^2$ a submersion for all $x \in \mathbb{T}^2.$ Indeed, denoting $y = \Phi^s_i(x)$, we have that
$$\left (\frac{\partial}{\partial t} \Phi_j^t \circ \Phi^s_i(x),  \frac{\partial}{\partial s} \Phi_j^t \circ \Phi^s_i(x)\right) = \left (D\Phi_j^t(y) F_j(y), D\Phi_j^t(y) F_i(y) \right).$$
 Proposition \ref{lem:submersivemap} implies that condition $(i)$ is satisfied. For the second assertion, we proceed as in the proof of Lemma \ref{lem:KLambda} $(iii).$ For all $\mu \in \C(\M)$ we have that
\begin{align*}
&\|\mu \Delta\|_\C(\M) \leq \|\mu\|_{\C(\M)}  \left[\int_\RR \max_{i = 1,2}  \|{\cal L}_{\Phi_i^t}\|_{C^{r-1}(M)} \alpha e^{-\alpha t} (1- \eta_n(t)) dt\right]^2\\
&\leq \left[\int_\RR C'  e^{-\beta t} (1- \eta_n(t)) dt\right]^2 \|\mu\|_{\C(\M)},
\end{align*}
for some constant $C', \beta > 0$ (by \tcr{(\ref{eq:boundalpha1}) and} (\ref{eq:boundnormL})).
For $n$ sufficiently large, the right-hand term can be made arbitrary small,  by monotone convergence.
\qed
\blem
\label{lem:Pnk}
We assume that $F_1,F_2$ are transverse at every point $p \in \mathbb{T}^2.$  Then for all $n \geq 2, k = 1, \ldots n-1,$ and $\eps > 0,$ $P_{n,k}$ can be decomposed into
$P_{n,k} = Q_{n,k} + \Delta_{n,k}$, where $Q_{n,k}, \Delta_{n,k}$ are Feller sub-Markov kernels and satisfy:
\bdes
\iti ${\cal M}(\M) Q_{n,k} \subset \C(\M)$;
\itii $\C(\M) \Delta_{n,k} \subset \C(\M)$;
\itiii for all $\mu \in \C(\M), \|\mu \Delta_{n,k}\|_{\C(\M)} \leq \eps \|\mu\|_{\C(\M)}.$
\edes
\elem
\prf With $Q, \Delta$ as in Lemma \ref{lem:KSK}, we set $$Q_{n,k} = (K \bar{S})^{k-1} Q A P^{n-k-1}, \quad \Delta_{n,k} = (K \bar{S})^{k-1}\Delta A P^{n-k-1}.$$ Then we have
$$P_{n,k}  = (K \bar{S})^{k-1} [KSK] A P^{n-k-1} = Q_{n,k} + \Delta_{n,k}.$$ Since ${\cal M}(\M)$ and $\C(\M)$ are invariant under the operators $K, A, \bar{S}, P,$ assertion $(i)$ and $(ii)$  follow directly from Lemma \ref{lem:KSK}.  We likewise have $\|\mu \Delta_{n,k}\|_{\C(\M)} \leq \eps \|{\cal K}\|_{C^{r-1}(\M)}^{n-2}  \|\mu\|_{\C(\M)}$, by Lemma \ref{lem:KSK}. Replacing $\eps$ by $\eps /  \|{\cal K}\|_{C^{r-1}(\M)}^{n-2}$, we obtain $(iii).$ \qed

\blem
\label{lem:Deltan}
We have the following:
\bdes
\iti $\C(\M) \Delta_{n,n} \subset \C(\M);$
\itii for all $\eps > 0,$ there exists $C<\infty$ such that for all $\mu \in \C(\M)$ and $n \geq 2$, we have
$$\|\mu \Delta_{n,n}\|_{\C(\M)} \leq C \left[\frac{\alpha - \min {(\alpha_{12}, \alpha_{21})}}{\alpha - \max_{i =1,2} \log \left({\cal R}({\cal L}_{\Phi_i^1}, C^{r-1}(M))\right)} \right]^n e^{n\eps} \|\mu\|_{\C(\M)}.$$
\edes
\elem
\prf We firstly observe that $K$ and $\bar{S}$ commute (\tcr{since the
rates are not position dependent}), so that $\Delta_{n,n} =\bar{S}^{n-1} K^n A.$  We therefore have that
\begin{align*}
&\|\mu \Delta_{n,n}\|_{\C(\M)} \leq \|\bar{S}\|^{n-1} \|{\cal K}^n\|_{C^{r-1}(\M)} \|A^t\| \|\mu\|_{\C(\M)}\\
&= \max(1- \frac{\alpha_{12}}{\alpha},1- \frac{\alpha_{21}}{\alpha})^{n-1} \|{\cal K}^n\|_{C^{r-1}(\M)}  \|A^t\| \|\mu\|_{\C(\M)}
\end{align*}
for all  $\mu \in \C(\M)$, whence the result follows from Lemma \ref{lem:KLambda} $(iii)$.
\qed

\bthm
\label{th:torus2}  Suppose that $F_1,F_2$ are $C^r, r \geq 1,$ transverse at every point $p \in \mathbb{T}^2,$ and that
 $$ \min {(\alpha_{12}, \alpha_{21})} >  \max_{i =1,2} \log \left({\cal R}({\cal L}_{\Phi_i^1}, C^{r-1}(M))\right).$$ Then every ergodic measure for $(Z_t)$ has a $C^{r-1}$ density with respect to $\m.$
\ethm
\prf Using the notation of the proceeding lemmas, we write
$P^n = Q_n + \Delta_n$, \tcb{where} $Q_n = \sum_{k = 0}^{n-1} Q_{n,k}$ and $\Delta_n = \sum_{k = 0}^{n} \Delta_{n,k}.$ Then $(Q_n, \Delta_n)$ satisfies the standing assumption, Assumption \ref{hyp1}, and for $n$ sufficiently large there exists $0 \leq \theta < 1$ such that $\|\mu \Delta_n\|_{\C(\M)} \leq \theta \|\mu\|_{\C(\M)}$ for all $\mu \in \C(\M).$ Theorem \ref{th:torus2} then follows from Theorem \ref{th:main2}
\qed
We then obtain Theorem \ref{th:torus} as a consequence of Theorem \ref{th:torus2} and the next proposition, Proposition \ref{prop:exptorus}, combined with the estimates given by Proposition \ref{prop:transfert}.

For a $C^1$ flow $\{\Phi^t\}$ we define the expansion rate and $k$-expansion volume rate of $\Phi$ to be the expansion rate and $k$-expansion volume rate of the time one map $\Phi^1$, which we denote by ${\cal E}(\Phi)$ and ${\cal EV}_k(\Phi)$ respectively.
\bprop
\label{prop:exptorus}
We let $F$ be a $C^2$ vector field on $\mathbb{T}^2$ with no equilibria (i.e $\mathsf{Eq}(F) = F^{-1}(0) = \emptyset$) and let $\{\Phi^t\}$ be the induced flow. Then
$${\cal E}(\Phi) =  \min_{\{\gamma \in \mathsf{Per}_{-}(F)\}} \Lambda_1(\gamma),$$ and
 $${\cal EV}_k(\Phi) = (k+1)\min_{\{\gamma \in \mathsf{Per}_{-}(F)\}} \Lambda_1(\gamma)$$ for all $k \geq 0,$ with the convention that the right-hand sides are $0$ whenever $\mathsf{Per}_{-}(F) = \emptyset.$
\eprop
\prf By  Propositions 14.2.2 and 14.2.4 in Katok and Hasselblat, \cite{KH95}, a fixed-point-free $C^2$ flow on $\mathbb{T}^2$ must enjoy one of the following two properties:
\bdes
\ita either all recurrent points are periodic;
\itb or there exists  a closed transversal and every orbit crosses this transversal. Furthermore, the return map to this transversal is a $C^2$ circle diffeomorphism  $f : S^1 \mapsto S^1$ which, by the Denjoy Theorem (\cite[Theorem 12.1.1]{KH95}), is topologically conjugate to an irrational rotation.
\edes
If $F$ has no periodic orbit then we are in case $(b)$.
We then have that ${\cal E}(\Phi) \leq 0$ by Remark \ref{rem:boundsE}. We now assume for contradiction that ${\cal E}(\Phi) < -\lambda < 0$. Then, by \cite[Corollary 2]{Sch97}, there exists two distinct points $x,y \in \mathbb{T}^2$ such that $\limsup_{t \rar \infty} \frac{\log(d(\Phi^t(x),\Phi^t(y))}{t} < -\lambda.$ This implies that the return map $f$ has two distinct points $\theta, \alpha \in S^1$ such that $d(f^n(\theta), f^n(\alpha) ) \rar 0$ as $n \rar 0$. However $f$ is topologically conjugate to a rotation and a rotation is an isometry, whence we obtain a contradiction.

If $F$ has periodic orbits, then we are in case $(a).$ We let $\mu$ be an ergodic probability measure for $\Phi^1.$ By the Poincaré recurrence theorem and Birkhoff's theorem, there exists a point $p$, recurrent for $\Phi^1,$ such that $$\frac{1}{n} \sum_{k = 1}^n \delta_{\Phi^k(p)} \Rightarrow \mu.$$ By $(a)$, $p$ is $T$-periodic for $\{\Phi^t\}$, for some $T > 0$. Thus, reasoning as in Example \ref{ex:OnS2}, either $\mu = \frac{1}{N} \sum_{i = 0}^{N-1} \delta_{\Phi^i(p)}$ for some $N\in \mathbb{N}$ (if $T$ is rational) or $\mu = \frac{1}{T} \int_0^T \delta_{\Phi^s(p)} ds$ (if $T$ is irrational). In both cases, $\Lambda_1(\mu)$ equals the Floquet exponent $\Lambda_1(\gamma)$ of the periodic orbit. The result then follows from Schreiber's theorem (equation (\ref{eq:Sch97})).
\qed

\brem
{\rm The fact that ${\cal E}(\Phi) = 0$ when $F$ has no periodic orbit  answers a question raised by Moe Hirsch in \cite{Hirsch94}. An affirmative  answer to this question is given in the introduction of Schreiber's paper \cite{Sch97}, but the proof and the assumptions are  not detailed in the paper. The result does actually directly  follows from Schreiber's results as shown above, at least for $C^2$ flows. The question is open for $C^1$ flows.}
\erem

\subsection{Smooth invariant distributions under fast switching}
\label{sec:fast}
We return here to the general model of a PDMP (as  described in the beginning of Section \ref{sec:pdmp}), but under the assumption that the rate matrix $(\alpha_{ij}(x))_{i,j \in E}$ is independent of $x$ and can be written as
\beq
\label{eq:simplerate}
\alpha_{ij}(x) =  \alpha a_{ij},
\eeq  where $\alpha > 0,$  $a_{ij} > 0$ for $i \neq j,$ and $a_{ii} = 0.$ The parameter $\alpha$ measures the rate of switching.

We shall prove here the following result.
\bthm
\mylabel{th:fastswitch}
Let $(Z_t)_{t \geq 0}$ be the PDMP corresponding to the characteristics $(\{F_i\}_{i \in E}, (\alpha_{ij})_{i,j \in E})$, where $\alpha_{ij}$ is given by (\ref{eq:simplerate}). Suppose that the $1$-Bracket condition holds at every point $x \in M.$ Then there exists $\alpha^* > 0$ such that, for all $\alpha \geq \alpha^*$, the ergodic measures of $(Z_t)$ (see Theorem \ref{th:lscpdmp1}) all have a $C^{r-1}$ density with respect to $\m.$
\ethm

A version of this result (under the assumption that there exists an accessible point), was established by the present authors in \cite{BT23}. However, the proof given here is simpler and provides a good illustration of our general method.

\subsubsection*{Proof of Theorem  \ref{th:fastswitch}}
Replacing $\alpha$ by $k \alpha$ and $a_{ij}$ by $a_{ij}/k,$ for $k$ sufficiently large, we can assume without loss of generality that $\sum_j a_{ij} <1.$ Set $A_{ij} := a_{ij}$ for $j \neq i$ and $A_{ii} := 1 - \sum_j a_{ij}.$

To highlight the influence of the switching rate parameter $\alpha$, we rewrite $K$ (as defined by (\ref{eq:K})) as $K_{\alpha}$ and $P$ (as defined by (\ref{eq:KLambda})) as
$$P_{\alpha} = K_{\alpha} A.$$ In light of Proposition \ref{prop:BLMZ24}, it suffices to consider invariant distributions of the operator $P_{\alpha}^n$ for some $n \geq 1.$

 For all $n \geq 2, \i = (i_1, \ldots, i_{n-1}) \in E^{n-1}$ and $i,j \in E,$ set
$$A[i, \i, j] =  A_{i i_1} A_{i_1 i_2}  \ldots A_{i_{n-2} i_{n-1}} A_{i_{n-1} j}$$ and
$$A[i,\i] = A_{i i_1} A_{i_1 i_2}  \ldots A_{i_{n-2} i_{n-1}} = \sum_j A[i,\i,j].$$
Let $h : (\RR_{+}^*)^{n} \mapsto [0,1]$ be a $C^{\infty}$ function and $\i \in E^{n-1}.$  Let  $P_{\alpha, \i, h}$ denote the sub-Markovian operator on $\M$
 defined by
$$P_{\alpha, \i, h} f(x,i) = \sum_{j \in E} A[i, \i, j]
 \int_{\RR^{n}_+} f(\Phi_{i_{n-1}}^{t_{n}/\alpha} \circ \cdots \circ \Phi_{i_1}^{t_2/\alpha} \circ \Phi_{i}^{t_1/\alpha}(x),j)   e^{-|t|}h(t) dt,$$ \tcb{where}
$|t| = t_1 + \ldots + t_{n}.$
 If $h \equiv 1,$  we write   $P_{\alpha, \i}$ for $P_{\alpha, \i, h.}$
Clearly we have that
$$P_{\alpha, \i} = P_{\alpha, \i, 1-h} + P_{\alpha, \i, h}$$
and
$$P_{\alpha, \i} f(x,i) = \EE_{x,i}(f(X_{n}, I_{n}) \Ind_{\left \{(I_1,\ldots, I_{n-1})  = \i \right \}}),$$ where $(X_n,I_n)$ is the discrete time Markov chain having $P_{\alpha}$ as transition kernel (see Remark \ref{rem:chainXI}).
In particular $$P_{\alpha}^{n}  = \sum_{\i \in E^{n-1}} P_{\alpha, \i}.$$

 Recall that $\|.\|_{C^k(M)}$ is a norm on $C^{k}(M)$ inducing the $C^k$ topology. For $\rho \in C^{k}(\M),$ define $\|\rho\|_{C^k(\M)}$ as $$\|\rho\|_{C^k(\M)}:= \sum_{i \in E} \|\rho_i\|_{C^k(M)}.$$
\blem We have the following:
\label{lem:fastswitch}
\bdes
\iti
 If $\mu \in {\cal M}_{ac}(\M)$ has density $\rho \in L^1(\m)$, then $\mu P_{\alpha,\i, h}$ has density ${\cal P}_{\alpha,\i, h}(\rho)$  given by
 $${\cal P}_{\alpha,\i, h}(\rho)(x,j) = \sum_i A[i, \i, j] {\cal L}_{\alpha,i,\i,h}(\rho_i)(x),$$
where
 $${\cal L}_{\alpha,i,\i,h}(\rho_i) := \int_{\RR^{n}_+} {\cal L}_{\Phi_{i_{n-1}}^{t_{n}/\alpha}} \circ \dots \circ {\cal L}_{\Phi_{i_1}^{t_2/\alpha}} \circ {\cal L}_{\Phi_{i}^{t_1/\alpha} }(\rho_i) e^{-|t|}h(t) dt.$$
 \itii If $\alpha > \max_{i \in E} \log \left({\cal R}({\cal L}_{\Phi_i^1}, C^{r-1}(M))\right)$, then ${\cal L}_{\alpha,i,\i,h}$ (respectively ${\cal P}_{\alpha,\i, h}$) is a bounded operator on $C^{r-1}(M)$ (respectively $C^{r-1}(\M)$).
 \itiii For all $\rho \in C^{r-1}(\M)$ and $\alpha > \max_{i \in E} \log \left({\cal R}({\cal L}_{\Phi_i^1}, C^{r-1}(M))\right)$,
 $$\|{\cal P}_{\alpha,\i, h} (\rho) \|_{C^{r-1}(\M)} \leq \epsilon_r(\alpha, h) \sum_{i \in E} A[i,\i] \|\rho_i\|_{C^{r-1}(M)}$$
 where
 $$\epsilon_r(\alpha,h):= \max_{i \in E, \i \in E^{n-1}} \|{\cal L}_{\alpha,i,\i,h} \|_{C^{r-1}(M)}.$$
Furthermore, for a convenient choice of norm $\|.\|_{C^{r-1}(M)},$
 $$\limsup_{\alpha \rar \infty} \epsilon_r(\alpha,h) \leq  \int_{\RR^n_+} e^{-|t|}h(t) dt.$$
 \edes
 \elem

\prf The proof of $(i)$ and $(ii)$ proceeds in the same manner as the proof of Lemma \ref{lem:KLambda} $(iii)$ (itself relying on Proposition \ref{lem:Lnu1}), so we refrain from repeating it for the sake of brevity.

$(iii).$ We have that
$$\|{\cal P}_{\alpha,\i,h}(\rho)\|_{C^{r-1}(\M)} = \sum_{j \in E} \|\sum_i A[i, \i, j] {\cal L}_{\alpha,i,\i,h}(\rho_i)\|_{C^{r-1}(M)}$$
$$\leq \sum_{j \in E} \sum_i A[i, \i, j] \epsilon_r(\alpha,h) \|\rho_i \|_{C^{r-1}(M)} =  \epsilon_r(\alpha,h) \sum_{i \in E} A[i,\i]  \|\rho_i \|_{C^{r-1}(M)}.$$
Let $\|. \|_{C^{r-1}(M)}$ be the $C_{r-1}$ norm induced by a finite atlas as in the proof of Proposition \ref{prop:transfert} (see equation (\ref{eq:cknormatlas})).

\medskip
\noindent
{\em Claim:} For all $j \in E, \limsup_{t \rar 0}  \| {\cal L}_{\Phi_{j}^{t}} \|_{C^{r-1}(M)} \leq 1.$

\medskip
\noindent
{\em Proof of the claim:} To shorten notation, set $\Phi^t = \Phi_j^t$ and  $F = F_j.$ Let
$$e_t(x) = \exp{[-\int_0^t \mathsf{div}(F)(\Phi^{-s}(x)) ds]}$$ and let $E_t, C_t$ be the operators defined by $$C_t(\rho)(x) = \rho(\Phi^t(x))$$ and $$E_t(\rho)(x) = e_t(x) \rho(x).$$ Thus, by formula (\ref{eq:LPhit}),
$${\cal L}_{\Phi^{t}} = E_t \circ C_t.$$ By the $C^r$ continuity of the map $(t,x) \mapsto \Phi^t(x)$ (see, for example, \cite[Chapter V, Corollary 4.1]{Hartman82}),
$\Phi^t \rar \Phi^0 = Id$ (the identity map), as $t \rar 0$ in the $C^r$ topology. Combined with Lemma \ref{lem:CDI} (ii) this implies that $\limsup_{t \rar 0} \|C_t\| \leq 1.$ This also implies that $e_t \rar 1,$ as $t \rar 0,$ in the $C^r$ topology, which combined with Lemma \ref{lem:CDI} (i), implies that  $\limsup_{t \rar 0} \|E_t\| \leq 1.$  This proves the claim.
\medskip

We let  $\eta(t_1, \ldots, t_n) =  \|{\cal L}_{\Phi_{i_{n-1}}^{t_n}}\|_{C^{r-1}(M)}  \dots \| {\cal L}_{\Phi_{i_0}^{t_1}} \|_{C^{r-1}(M)},$ where here $i_0$ stands for  $i.$ It follows from the claim that $\limsup_{t \rar 0_{\RR^n}} \eta(t) \leq 1.$
Therefore for all $\eps > 0$, there exists some $\eps > 0$ such that $\eta(t) \leq 1 + \eps$ for all $t \in \RR^n_+$ such that $|t| \leq \delta.$ Thus we have
\begin{align*}
&\|{\cal L}_{\alpha,\i,h} \|_{C^{r-1}(M)} \leq \int_{\RR^n_+} \eta(t/\alpha)  e^{-|t|} h(t)  dt\\
&\leq (1+\eps) \int_{\RR^n_+} e^{-|t|} h(t) \Ind_{|t| \leq \alpha \delta} dt   +  \int_{\RR^n_+ } \eta(t/\alpha) h(t)  e^{-|t|} \Ind_{|t| \geq \alpha \delta} dt.
\end{align*}
When $\alpha \rar \infty,$ the first term on the right goes to $1+\eps$ while the second term goes to $0$. This follows from the fact that $\eta(t) \leq C' e^{\beta |t|}$ for some $\beta > 0$ and $C'<\infty$, by equation (\ref{eq:boundnormL}). This concludes the proof. \qed
%%%%%%%%%%%%
\bprop
\label{prop:fastswitch} We suppose that there exist $n \geq 2, \i = (i_1, \ldots, i_{n-1}) \in E^{n-1}$ and $U \subset (\RR_{+}^*)^{n-1}$ a nonempty open set such that:
\bdes
\iti $\frac{1}{\alpha} U \subset U$ for all $\alpha \geq 1;$
\itii for all $x \in M,$ the map  $(t_2, \ldots, t_n) \rar  \Phi^{t_n}_{i_{n-1}} \circ \dots \circ \Phi^{t_2}_{i_1}(x)$ is a submersion on $U.$
\edes
Then, there exists $\alpha^* \geq 1$ such that $\mathsf{Inv}(P_{\alpha}) \subset {\cal M}_{ac}^{r-1}(M)$ for all $\alpha \geq \alpha^*$.
\eprop
\prf We let  $h : (\RR_{+}^*)^n \rar [0,1]$ be a $C^{\infty}$ non-identically zero function with compact support in $R_+^* \times U.$ We set $Q_{\alpha} = P_{\alpha, \i, h}$ and $\Delta_{\alpha} = P_{\alpha, \i, 1-h} + \sum_{\mathbf{j} \in E^{n-1} \setminus{ \{\i}\}} P_{\alpha, \mathbf{j}}$,
and take ${\cal C}(M) = {\cal M}_{ac}^{r-1}(M).$

The conditions $(i)$ and $(ii)$ imply that the map $$(t_1,t_2 \ldots, t_n) \mapsto \Phi^{t_n/\alpha}_{i_{n-1}} \circ \dots \circ \Phi^{t_2/\alpha}_{i_1} \circ \Phi^{t_1/\alpha}_{i} (x)$$ is a submersion on $\RR_{+}^* \times U,$ for all $\alpha > 0.$ Thus, by Proposition \ref{lem:submersivemap}, $${\cal M}(M) Q_{\alpha} \subset {\cal C}(M)$$ for all $\alpha \geq 1.$  By Lemma \ref{lem:fastswitch}, for $\alpha$ sufficiently large
and for all $\rho \in C^{r-1}(\M),$
$$\|{\cal P}_{\alpha, \i, 1-h}(\rho) + \sum_{\mathbf{j} \in E^{n-1} \setminus{ \{\i}\}} {\cal P}_{\alpha, \mathbf{j}}(\rho)\|_{C^{r-1}(\M)}$$
$$\leq \epsilon_r(\alpha,1-h) \sum_{i \in E} A[i,\i] \|\rho_i\|_{C^{r-1}(M)} + \epsilon_r(\alpha,1) \sum_{i \in E} (1- A[i,\i])  \|\rho_i\|_{C^{r-1}(M)}.$$
For all  $\epsilon > 0$ there exists $\alpha^*$ such that  $\epsilon_r(\alpha,1-h) \leq 1 - \int_{(\RR_{+}^*)^n} e^{- |t|} h(t)  dt + \epsilon$ and
$\epsilon_r(\alpha,1) \leq 1 + \epsilon$ for all $\alpha \geq \alpha^*.$ It then follows that, for all $\alpha \geq \alpha^*,$
$$
\|\Delta_{\alpha}\|_{{\cal C}(M)} :=  \|{\cal P}_{\alpha, \i, 1-h} + \sum_{\mathbf{j} \in E^{n-1} \setminus{ \{\i}\}} {\cal P}_{\alpha, \mathbf{j}}\|_{C^{r-1}(\M)}$$
$$ \leq [1- \min_{i \in E} A[i,\i] \int_{(\RR_{+}^*)^n} e^{- |t|} h(t)  dt] + O(\epsilon).$$ For $\epsilon$ sufficiently small, this latter quantity is $< 1$ and the proposition then follows from Theorem \ref{th:main2}.
\qed

By \cite[Proposition 5.1]{BT23}, the $1$-Bracket condition implies that the assumptions of Proposition \ref{prop:fastswitch} are satisfied. This concludes the proof of Theorem \ref{th:fastswitch}.\qed
\medskip
\subsection{PDMPs on noncompact manifolds}
Suppose that the vector fields $F_i$ are defined on a (possibly) noncompact $d$-dimensional manifold $W$ (typically $\RR^d$), and that there exists a compact connected  $d$-dimensional submanifold $M \subset W$ with nonempty boundary $\partial M$ such that for each $x \in \partial M$ and $i \in E,$ $F_i(x)$ points inward $M.$ Then  all the preceding results  remain valid  for the PDMP living on $M.$
\bex
{\rm
This simple example generalizes Example 4.7 given by Malrieu (\cite{Malrieu}) and provides a partial answer to his {\em Open Question 4}.

 Let $d \geq 2.$ Let $A$ be a $d \times d$ real matrix which is not a dilation, whose eigenvalues have all negative real parts. Let $H \subset \RR^d$ be a $d-1$ dimensional vector space such that $A H \neq H.$ Let $p_1, \ldots p_{d-1}$ be a basis of $H,$ and set $p_d := 0.$ Define affine vector fields $F_1, \ldots, F_d$ on $\RR^d$ by $$F_i(x) = A(x-p_i).$$
Because eigenvalues of $A$ have negative real parts, there exist $r > 0$ and $\tau > 0$ such that
$$\|e^{tA} x\| \leq e^{-r t} \|x\|$$ for all $t  \geq \tau,$ where $\|x\| = \sqrt{\langle x, x \rangle}$ is the standard Euclidean norm of $x \in \RR^d.$  Let  $$\left( x, y \right) = \int_0^{\tau} e^{2 r s} \langle e^{sA} x, e^{sA} y \rangle ds$$ and $V(x) = \sqrt{\left( x, x \right)}.$ Then,  $V$ is an {\em adapted Euclidean} norm  on $\RR^d$ in the sense  that $$V(e^{tA} x) \leq e^{-r t} V(x)$$ {\bf for all} $t \geq 0$ and $x \in \RR^d$ (see for instance the proof   Theorem 5.1 in  \cite{Robinson}). Thus, for all $x \neq 0,$
$$\lim_{t \rar 0} \frac{V(e^{tA} x) - V(x)}{t} = \left( F_d(x), \nabla V(x) \right) \leq - r V(x)$$  and, for all $i = 1, \ldots d-1,$ $$\left( F_i(x), \nabla V(x) \right) \leq - r V(x) - \left( A p_i, \nabla V(x) \right) $$
$$ \leq - r V(x) + V(Ap_i).$$
Fix $R > \max_{i = 1, \ldots d-1} \frac{V(Ap_i)}{r}$ and let $$M= \{x \in \RR^d:   \: V(x) \leq R\}.$$ Then $M$ is a compact submanifold of $\RR^d$ with boundary $\partial M = V^{-1}(R),$ and each $F_i$ points inward $M$ at $\partial M.$

We claim that the $1$-Bracket condition holds true at every point $x \in M.$ Indeed,  elementary computations show that
 $$[F_i, F_d](x) = A^2 p_i,$$
$$\mathsf{det}([F_1,F_d](x),\ldots, [F_{d-1}, F_d](x), F_d(x)) = \mathsf{det}(A) \mathsf{det}(Ap_1, \ldots, Ap_{d-1}, x),$$
and, for all $k = 1, \ldots, d-1,$
$$\mathsf{det}([F_1,F_d](x),\ldots, [F_{d-1}, F_d](x), F_k(x))$$
  $$=\mathsf{det}(A) (\mathsf{det}(Ap_1, \ldots, Ap_{d-1}, x) - \mathsf{det}(Ap_1, \ldots, Ap_{d-1}, p_k))$$
If the first determinant is nonzero, the condition holds. If it is zero, pick $k = 1, \ldots d-1$ such that $p_k \not \in A H$ (recall that $H \neq AH$).
For such a $k$ the second determinant is nonzero.

 Consider now the PDMP on $\M = M \times E$ with $E = \{1, \ldots, d\}$ having characteristics
$(\{F_i\}_{i \in E},   (\alpha a_{ij})_{ij \in E})$ with $\alpha > 0$ and $a_{ij} > 0$ for all $i \neq j.$ One has the following properties:

\bdes
\iti The PDMP $(Z_t)$ has a unique invariant probability $\Pi$ absolutely continuous with respect to the Lebesgue measure whose density $\rho$ is lower semi continuous  with respect to Lebesgue. This follows from Theorem \ref{th:lscpdmp0} because the origin (or any point $p_i$) is accessible and satisfies the weak bracket condition.

\itii  For $\alpha$ sufficiently large,  $\rho$ is $C^k$ by Theorem \ref{th:fastswitch}.

 \itiii Furthermore it can be shown (see Theorem 2.13  in \cite{BT23}) that
 $$\rho_i(p_i) = \infty$$
for $$\alpha \sum_{j \neq i} a_{ij} \leq - \mathsf{Tr}(A).$$
 \edes
Observe that if  $H  =  AH,$  there is still a unique invariant measure (because the flows induced by the $F_i$ contract distances) which is supported by $H,$ hence necessarily singular with respect to Lebesgue.
}
\eex
{\textbf{Acknowledgement:}}  The work of MB, and partially that of OT, was funded by the grant 200020-219913 from the Swiss National Foundation. The work of OT was also partially funded by the EPSRC MathRad programme grant EP/W026899/.
We thank two anonymous referees for their valuable comments and suggestions.

Conflict of Interest : None.
\bibliographystyle{amsplain}

\bibliography{regulardensity.bib}

\providecommand{\bysame}{\leavevmode\hbox to3em{\hrulefill}\thinspace}
\providecommand{\MR}{\relax\ifhmode\unskip\space\fi MR }
\providecommand{\MRhref}[2]{%
  \href{http://www.ams.org/mathscinet-getitem?mr=#1}{#2}
}
\providecommand{\href}[2]{#2}
\begin{thebibliography}{10}

\bibitem{AN78}
K.~B. Athreya and P.~Ney, \emph{A new approach to the limit theory of recurrent
  markov chains}, Transactions of the American Mathematical Society
  \textbf{245} (1978), 493--501.

\bibitem{bakhtin&hurt&matt}
Y.~Bakhtin, T.~Hurth, and J.C. Mattingly, \emph{Regularity of invariant
  densities for 1d-systems with random switching}, Nonlinearity \textbf{28}
  (2015), 3755--3787.

\bibitem{bakhtin&hurt}
Yuri. Bakhtin and Tobias. Hurth, \emph{Invariant densities for dynamical
  systems with random switching}, Nonlinearity (2012), no.~10, 2937--2952.

\bibitem{BHLM18}
Yuri Bakhtin, Tobias Hurth, Sean~D. Lawley, and Jonathan~C. Mattingly,
  \emph{Smooth invariant densities for random switching on the torus},
  Nonlinearity \textbf{31} (2018), no.~4, 1331--1350. \MR{3816636}

\bibitem{BHLM21}
\bysame, \emph{Singularities of invariant densities for random switching
  between two linear {ODE}s in 2{D}}, SIAM J. Appl. Dyn. Syst. \textbf{20}
  (2021), no.~4, 1917--1958. \MR{4322095}

\bibitem{Baladi2000}
Viviane Baladi, \emph{Positive transfer operators and decay of correlations},
  Advanced Series in Nonlinear Dynamics, vol.~16, World Scientific Publishing
  Co., Inc., River Edge, NJ, 2000. \MR{1793194}

\bibitem{BHS18}
M.~Bena\"im, T.~Hurth, and E.~Strickler, \emph{A user-friendly condition for
  exponential ergodicity in randomly switched environments}, Electronic
  Communications in Probability \textbf{23} (2018), no.~44, 1--12.

\bibitem{BLMZ_2015}
M.~Bena\"{\i}m, S.~Le~Borgne, F.~Malrieu, and P.-A. Zitt, \emph{Qualitative
  properties of certain piecewise deterministic {M}arkov processes}, Ann. Inst.
  Henri Poincar\'{e} Probab. Stat. \textbf{51} (2015), no.~3, 1040--1075.

\bibitem{BLMZ_2019}
\bysame, \emph{Erratum: Qualitative properties of certain piecewise
  deterministic {M}arkov processes}, Ann. Inst. Henri Poincar\'{e} Probab.
  Stat. \textbf{55} (2019), no.~4.

\bibitem{BH95}
Michel Bena\"{\i}m and Morris~W. Hirsch, \emph{Chain recurrence in surface
  flows}, Discrete Contin. Dynam. Systems \textbf{1} (1995), no.~1, 1--16.
  \MR{1355862}

\bibitem{BenaimHurth}
Michel Bena\"{\i}m and Tobias Hurth, \emph{Markov chains on metric spaces, a
  short course}, Universitext, vol.~99, Springer, Cham, 2022.

\bibitem{BCL17}
M~Benaïm, F~Colonius, and R~Lettau, \emph{Supports of invariant measures for
  piecewise deterministic {M}arkov processes}, Nonlinearity \textbf{30} (2017),
  no.~9, 3400.

\bibitem{BT23}
Michel Benaïm and Oliver Tough, \emph{Regularity of the stationary density for
  systems with fast random switching}, 2023.

\bibitem{campbell}
J.~Campbell and Latushkin Y., \emph{Sharp estimates in {R}uelle theorems for
  matrix transfer operators}, Comm. Math. Phys. \textbf{185} (1997), 379--396.

\bibitem{Cloez15}
Bertrand Cloez and Martin Hairer, \emph{Exponential ergodicity for {M}arkov
  processes with random switching}, Bernoulli \textbf{21} (2015), no.~1,
  505--536. \MR{3322329}

\bibitem{MR1283589}
M.~H.~A. Davis, \emph{Markov models and optimization}, Monographs on Statistics
  and Applied Probability, vol.~49, Chapman \& Hall, London, 1993.

\bibitem{diestel}
J.~Diestel and J.~J. Uhl, Jr., \emph{Vector measures}, Mathematical Surveys,
  vol. No. 15, American Mathematical Society, Providence, RI, 1977, With a
  foreword by B. J. Pettis. \MR{453964}

\bibitem{duf00}
Marie Duflo, \emph{Random iterative models}, Applications of Mathematics (New
  York), vol.~34, Springer-Verlag, Berlin, 1997, Translated from the 1990
  French original by Stephen S. Wilson and revised by the author. \MR{1485774}

\bibitem{Hartman82}
Philip Hartman, \emph{Ordinary differential equations}, second ed.,
  Birkh\"{a}user, Boston, Mass., 1982. \MR{658490}

\bibitem{Hirsch76}
Morris~W. Hirsch, \emph{Differential topology}, Graduate Texts in Mathematics,
  No. 33, Springer-Verlag, New York-Heidelberg, 1976. \MR{0448362}

\bibitem{Hirsch94}
\bysame, \emph{Asymptotic phase, shadowing and reaction-diffusion systems},
  Differential equations, dynamical systems, and control science, Lecture Notes
  in Pure and Appl. Math., vol. 152, Dekker, New York, 1994, pp.~87--99.
  \MR{1243195}

\bibitem{KH95}
Anatole Katok and Boris Hasselblatt, \emph{Introduction to the modern theory of
  dynamical systems}, Encyclopedia of Mathematics and its Applications,
  vol.~54, Cambridge University Press, Cambridge, 1995, With a supplementary
  chapter by Katok and Leonardo Mendoza. \MR{1326374}

\bibitem{Loch18}
E.~L\"{o}cherbach, \emph{Absolute continuity of the invariant measure in
  piecewise deterministic {M}arkov processes having degenerate jumps},
  Stochastic Process. Appl. \textbf{128} (2018), no.~6, 1797--1829.
  \MR{3797644}

\bibitem{Mane87}
Ricardo Ma\~{n}\'{e}, \emph{Ergodic theory and differentiable dynamics},
  Ergebnisse der Mathematik und ihrer Grenzgebiete (3) [Results in Mathematics
  and Related Areas (3)], vol.~8, Springer-Verlag, Berlin, 1987, Translated
  from the Portuguese by Silvio Levy. \MR{889254}

\bibitem{Malrieu}
Florent Malrieu, \emph{Some simple but challenging {Markov} processes}, Annales
  de la Facult\'e des sciences de Toulouse : Math\'ematiques \textbf{Ser. 6,
  24} (2015), no.~4, 857--883 (en).

\bibitem{MT93}
Sean Meyn and Richard~L. Tweedie, \emph{Markov chains and stochastic
  stability}, second ed., Cambridge University Press, Cambridge, 2009, With a
  prologue by Peter W. Glynn. \MR{2509253}

\bibitem{Os68}
Valery Oseledets, \emph{A multiplicative ergodic theorem. characteristic
  ljapunov, exponents of dynamical systems}, Trudy Moskov. Mat. Obsch
  \textbf{19} (1968), 179--210.

\bibitem{Robinson}
Clark Robinson, \emph{Dynamical systems}, second ed., Studies in Advanced
  Mathematics, CRC Press, Boca Raton, FL, 1999, Stability, symbolic dynamics,
  and chaos. \MR{1792240}

\bibitem{Ruelle89}
David Ruelle, \emph{The thermodynamic formalism for expanding maps}, Comm.
  Math. Phys. \textbf{125} (1989), no.~2, 239--262. \MR{1016871}

\bibitem{Sch97}
Sebastian~J. Schreiber, \emph{Expansion rates and {L}yapunov exponents},
  Discrete Contin. Dynam. Systems \textbf{3} (1997), no.~3, 433--438.
  \MR{1444204}

\bibitem{Sch98}
\bysame, \emph{On growth rates of subadditive functions for semiflows}, J.
  Differential Equations \textbf{148} (1998), no.~2, 334--350. \MR{1643183}

\end{thebibliography}
\end{document}